\documentclass[10pt,journal]{IEEEtran}

\usepackage[utf8]{inputenc}
\usepackage{amsmath}
\usepackage{amsthm}
\usepackage{amssymb}
\usepackage{graphicx}
\usepackage{psfrag}
\usepackage{cite}
\usepackage{url}
\usepackage{enumerate}
\usepackage[usenames]{color}

 \theoremstyle{plain}
 \newtheorem{thm}{Theorem}
 \newtheorem{lem}[thm]{Lemma}
 \newtheorem{prop}[thm]{Proposition}
 \newtheorem{cor}[thm]{Corollary}
 \theoremstyle{definition}
 \newtheorem{defn}{Definition}
 \newtheorem{rem}{Remark}

 \newtheorem{assm}{Assumption}

\newcommand{\bb}[0]{\begin{bmatrix}}
\newcommand{\eb}[0]{\end{bmatrix}}
\newcommand{\be}[0]{\begin{equation}}
\newcommand{\ee}[0]{\end{equation}}
\newcommand{\ben}[0]{\begin{equation*}}
\newcommand{\een}[0]{\end{equation*}}

\newcommand{\norms}[1]{\lVert#1\rVert}
\newcommand{\normB}[0]{\norms}
\let\norm=\normB

\newcommand{\abs}[1]{\left\lvert#1\right\rvert}

\renewcommand{\Re}[0]{\mathbb R}

\newcommand{\norml}[2]{\norm{#1}_{L_{#2}}}

\DeclareMathOperator*{\argmin}{arg\,min}

\renewcommand{\IEEEQED}{\IEEEQEDopen}

\title{Adaptive Systems with Closed--loop Reference Models: Stability, Robustness and Transient Performance}

\author{Travis~E.~Gibson,
Anuradha~M.~Annaswamy and
        Eugene~Lavretsky
\thanks{T.~E. Gibson and A. M. Annaswamy are with the Department
of Mechanical Engineering, Massaschusetts Institute of Technology, Cambridge,
MA, 02139 e-mail: ({tgibson@mit.edu}).}
\thanks{E. Lavretsky is with The Boeing Company, Huntington Beach, CA 92648}
\thanks{This work was supported in part by the Boeing Strategic Research Initiative and in part by the National Institute for Aerospace Activity 2621.}}

\begin{document}

\maketitle

\begin{abstract}
This paper explores the properties of adaptive systems with closed--loop reference models. Using additional design freedom available in closed--loop reference models, we design new adaptive controllers that are (a) stable, and (b) have improved transient properties. Numerical studies that complement theoretical derivations are also reported.
\end{abstract}

\section{Introduction}
The central element of any adaptive systems is online parameter adjustment. This is usually accomplished by having a plant, determined by a dynamic model, along with a controller with adaptive parameters designed to compensate for the plant's actions, follow a reference model. The resulting error between the reference model and the plant is used to adjust the adaptive parameter. 


By definition, open--loop reference models are independent of the system dynamics. Such reference models have been the backbone of adaptive control for the past four decades \cite{annbook,ioabook} where modifications to the adaptive control law were first added for stability in the presence of bounded disturbances \cite{kre82,ioa83book,narendra86} and semi--global stability in the presence of unmodeled dynamics\cite{tsa86,nar87TAC}.  We denote the underlying open--loop systems in all these cases as {\em Open-loop Reference Model} (ORM)--adaptive systems.

Earlier developments of adaptive systems included explorations of various kinds of reference models. The overall goal behind the selection of a reference model is that the corresponding {\em tracking error} must asymptotically decay in the absence of parametric uncertainties in the plant being controlled. In order to accomplish this goal, modifications of the open--loop reference models were explored \cite{lan74,lan79}. Some of these modifications retained stability properties and were otherwise indistinguishable from ORM--adaptive systems and as a result, not pursued. Others could not be shown to be stable and were therefore dropped. Recently, a class of {\em Closed--loop Reference Models} (CRM) have been proposed for control of plants with unknown parameters whose states are accessible (see for example \cite{lav10aiaa,lav12tac,yale10,ste10,ste11}) all of which are guaranteed to be stable and additionally portray improved transient performance.

Transient  performance has been directly addressed in  \cite{krs93} and more recently in\cite{lav12tac,hov10,ste10,ste11} . The results in \cite{lav12tac} discussed the tracking error, but focused the attention mainly on the initial interval where the CRM-adaptive system exhibits fast time-scales. Results in \cite{ste10,ste11} focus on deriving a damping ratio and natural frequency for adaptive systems with CRM. However, assumptions are made that the initial state error is zero and that the closed-loop system state is independent of the feedback gain in the reference model, both of which may not hold in general. The results in \cite{hov10} too assume that initial state errors are zero. And in addition, the bounds derived in \cite{hov10} are based upon $\mathcal L_\infty$ norms, which do not capture the transient properties of adaptive systems. The results in \cite{krs93} pertain to transient properties of adaptive systems, and quantify them using an $\mathcal L_2$ norm. The adaptive systems in question however are indirect, and do not pertain to CRMs. With the exception of \cite{ste10,ste11}, none of the others have focused on derivatives of signals in the adaptive system, which is another measure of transient performance.


Our focus in this is paper on CRM-based adaptive systems. Similar to \cite{lav12tac,ste10,ste11} we demonstrate their stability properties. Unlike these papers, we discuss transient properties of these adaptive systems using an ${\mathcal L_2}$ norm of error signals and derivatives of key signals such as adaptive parameters and the control input. These metrics are used to compare the CRM adaptive systems with their ORM counterparts

%

Another class of adaptive systems that have been explored in the past where a noticeable improvement in transient performance was obtained is in the context of {\em Combined/composite direct and indirect Model Reference Adaptive Control} (CMRAC) \cite{duarte1989combined,slotine1989composite}. While the results of these papers established stability of combined schemes, no rigorous guarantees of improved transient performance were provided, and have remained a conjecture \cite{eugeneTAC09}.  We focus on this class of adaptive systems as well in this paper and introduce CRMs into the picture. The resulting {\em CMRAC with Closed-loop reference models} (CMRAC-C) are shown to be stable, enable the feedback of noise-free state estimates while guaranteeing stability, and most importantly are shown to have guaranteed transient properties.


The main contributions of the paper are (i) direct adaptive control structures with guaranteed transient performance, (ii) combined direct and indirect adaptive controllers with guaranteed transient performance, and (iii) the development of adaptive systems that allow feedback from noise free regressors. These are realized by using the extra degree of freedom available in the CRM in terms of a feedback gain, and by exploiting exponential convergence properties of the CRM--adaptive system. The latter is made possible by introducing a projection algorithm with a known upper bound on the unknown parameters. 

Previous work by the same authors on this subject can be found in \cite{gib12,gib13acc1,gib13acc2}. The findings presented in   \cite{gib12} are preliminary and illustrate the waterbed affect through Euclidean bounds of the rate of control input at specific times of interest. The results presented in \cite{gib13acc1,gib13acc2} are a condensed version of this article and do not contain all of the discussions or proofs necessary to have a clear self contained presentation of this material. While much of what is presented in the paper is restricted to plants with state-variables accessible, the same idea can be extended to adaptive control using output feedback and is the topic of current investigation \cite{gib13ecc}.

The results in this paper are organized as follows: Section II introduces the basic structure of CRM adaptive control as well as the Projection Operator. Section III investigates the transient response of CRM. Section IV investigates the robustness properties of CRM adaptive control. Section V contains the stability analysis of CMRAC--C. Section VI analyses the transient performance of CMRAC--C. Section VII contains the analysis of CMRAC--CO and its robustness properties in regard to measurement noise. Section VIII compares CRM, CMRAC-C and CMRAC--CO structures. Section IX contains our concluding remarks.

\section{The CRM--Adaptive System}\label{sec:l}

In this section, we describe the CRM--adaptive system, and establish its stability and convergence properties in the absence of any perturbations other than parametric uncertainties. We first describe the CRM--adaptive system and prove its closed--loop stability. After some preliminaries on matrix bounds, we introduce a projection algorithm in the adaptive law. This is used to derive exponentially converging bounds on the key variables in the CRM--adaptive system. 

Consider the linear system dynamics with scalar input
\be\label{s1:sys}
\dot x(t) = A_p x(t) + b u(t)
\ee
where ${x(t) \in \Re^n}$ is the state vector, ${u(t)\in \Re}$ is the control input, ${A_p\in \Re^{n \times n}}$ is unknown and ${b\in \Re^n}$ is known. Our goal is to design the control input such that ${x(t)}$ follows the reference model state ${x_m(t)\in \Re^n}$ defined by  the following dynamics
\be\label{s1:ref}
\dot x_m(t) = A_m x_m(t) + b r(t) - L(x(t)-x_m(t))
\ee
where ${A_m \in \Re^{n\times n }}$ is Hurwtitz and ${r(t)\in \Re}$ is a bounded possible time varying reference command. ${L\in \Re^{n\times n}}$ is denoted as the {\em Luenberger--gain}, and is chosen such that 
\be\label{s1:ambar}
\bar A_m \triangleq A_m+L
\ee
is Hurwitz. Equation \eqref{s1:ref} is referred to as a CRM, and when ${L=0}$ the classical ORM is recovered.

\begin{assm}\label{asm:m}
A parameter vector ${\theta^*\in\Re^n}$ exists that satisfies the {\em matching condition}
\be\label{s1:match}
A_m=A_p+b\theta^{*T}.
\ee
\end{assm}

The control input is chosen in the form
\be\label{s1:u}
u(t)= \theta^T(t) x(t) + r(t)
\ee
where $\theta(t) \in \Re^n$ is the adaptive control gain with the update law
\be\label{s1:law}
\dot \theta(t) = -\Gamma x(t)e^T(t) Pb
\ee
with ${\Gamma = \Gamma^T>0}$, ${e(t)= x(t)-x_m(t)}$ is the model following error and $P=P^T>0$ is the solution to the algebraic Lyapunov equation
\be
\label{s1:lyap}\bar A_m^TP+P \bar A_m = -I_{n\times n}. 
\ee
The underlying error model in this case is given by
\be\dot e(t)= \bar A_m e(t) + b \tilde\theta(t)x(t)   \label{s1:et}\ee
where ${\tilde\theta(t)=\theta(t)-\theta^*}$ is the parameter error.

\begin{thm}\label{thm:simple}
The closed-loop adaptive system with \eqref{s1:sys}, \eqref{s1:ref}, \eqref{s1:u} and \eqref{s1:law} is globally stable with $e(t)$ tending to zero asymptotically, under the matching condition in \eqref{s1:match}.
\end{thm}
\begin{proof} It is straight forward to show using \eqref{s1:law} and \eqref{s1:et} that
\be\label{can}
V(e,\tilde\theta) = e^TPe + \tilde \theta^T \Gamma^{-1} \tilde \theta
\ee
is a Lyapunov function. Since $e$ is bounded, the structure of \eqref{s1:ref} implies that $x_m$ is bounded. $x$ in turn and $u$ are bounded. Barbalat  lemma ensures asymptotic convergence of $e(t)$ to zero. \end{proof}

\begin{cor}\label{Tdef}
For all $\epsilon>0$ there exists $T(\epsilon,L)>0$ such that $t\geq T(\epsilon,L)$ implies $\norm{e(t)}\leq \epsilon$.
\end{cor}
\noindent The above corollary seems redundant. It is however used later in the transient performance.
%
%

The overall CRM--adaptive system is defined by \eqref{s1:sys}, \eqref{s1:ref}, \eqref{s1:u}, and \eqref{s1:law}.  The standard open-loop reference model is given by 
\be\dot x_m^o(t) = A_m x_m^o(t) + b r(t) \label{eo}\ee
with the corresponding tracking error 
\be e^o(t) =  x(t) - x_m^o(t) \label{eodef}. \ee
One can in fact view the error $e^o$ as the {\em true tracking error} and $e$ as a pseudo--tracking error. The question that arises is whether the convergence properties that are assured in an ORM--adaptive system, of $e^o(t)$ tending to zero is guaranteed in a CRM--adaptive system as well. This is addressed in the following corollary:

\begin{cor}
The state vector $x(t)$ converges to $x_m^o(t)$ as $t\rightarrow \infty$.
\end{cor}
\begin{proof}
From Theorem \ref{thm:simple} we can conclude that $e(t)\rightarrow0$ asymptotically. Thus we can conclude that  $x_m(t) \rightarrow x_m^o(t)$ as $e(t)\rightarrow 0$, implying that $e^o\rightarrow 0$, thus $x(t)\rightarrow x_m^o(t)$ as $t\rightarrow \infty$. 
\end{proof}

\begin{rem}
The choice of the CRM as in (2) essentially makes the reference model nonlinear, as $x$ depends on $\theta$ which in turn depends on $x_m$ in  a highly nonlinear manner. In general, the tracking problem in a ORM--adaptive system can be viewed as one where an overall nonlinear time-varying system is to be designed such that its output tracks that of a linear time-invariant system. The CRM-adaptive system is one where the overall nonlinear system is instead required to follow a nonlinear (reference) model. This nonlinear model, however, is chosen such that it asymptotically approaches the original linear reference model in the classical case, and hence the CRM-adaptive system retains all the desired characteristics of the ORM-adaptive system. As we will show in Section \ref{section:perf}, the CRM-adaptive system has an additional desirable property, of quantifiable transient properties, which the ORM-adaptive system may not necessarily possess. We will also show in this section that this is made possible by virtue of the additional degree of freedom available to the adaptive system in the form of the feedback gain in the CRM.
\end{rem}

%

\subsection{Preliminaries}
All norms unless otherwise noted are the Euclidean--norm and the induced Euclidean--norm. The variable $t\in \Re_+$ denotes time throughout and for a differentiable function $x(t)$, $\frac{d}{dt} x(t)$ is equivalent to $\dot x(t)$. Parameters explicit time dependence $(t)$ is used upon introduction and then omitted thereafter except for emphasis. The other norms used in this work are the $\mathcal L_2$ and truncated $\mathcal L_2$ norm defined below. Given a vector $\nu \in \Re^n$ and finite $p \in\mathbb N_{>0}$, $
\norm {\nu(t)}_{L_p}    \triangleq  \left(  \int_0^\infty \norm {\nu(s)} ^p ds\right)^{1/p}$, $\norm {\nu(t)}_{L_p,\tau}  \triangleq   \left( \int_0^\tau \norm {\nu(s)} ^p ds\right)^{1/p}$ and the infinity norm is then defined as $\norml{\nu(t)}{\infty} \triangleq \sup \norm{\nu(t)}$.

\begin{defn} Given a Hurwtiz matrix $A_m\in \Re^{n\times n}$
\be
\begin{split}
\label{s1:ambound}
 \sigma& \triangleq - \max_i  \left( \text{real} (\lambda_i(A_m))  \right) \\
  s& \triangleq - \min_i  \left( \lambda_i\left(A_m+A_m^T\right)/2  \right) \\
 	  a& \triangleq \norm{A_m}. 
\end{split}
\ee 
 For ease of exposition, throughout the paper, we choose $L$ in  \eqref{s1:ref} and $\Gamma$ in \eqref{s2:adaplaw} as follows:
\begin{align}\label{s1:L} L &\triangleq  -\ell I_{n \times n} \\
\label{s1:G}\Gamma & \triangleq \gamma I_{n\times n}.
\end{align}
\end{defn}

\begin{lem}\label{sgs}
The constants $\sigma$ and $s$ are strictly positive and satisfy
\ben
s\geq\sigma>0.
\een
\end{lem}
\begin{proof}
 $A_m$ is Hurwitz and therefore $\sigma>0$. It is not necessary however that the sum $A_m+A_m^T$ is Hurwitz. The trace operator is denoted as $\text{Tr} (\cdot)$ and is a linear operator. Recalling that a matrix and its transpose have the same trace we can conclude that $\Sigma _i \lambda_i(A_m) =\Sigma _i \lambda_i(A_m+A_m^T)/2$. Finally we have that
 \ben
 s \geq -\frac 1 {2n} \sum_{i=1}^n \lambda_i\left(A_m+A_m^T\right) = -\frac 1 n \sum_{i=1}^n \lambda_i(A_m) \geq \sigma>0.  \ \IEEEQED
 \een
 \let\IEEEQED\relax 
\end{proof}

\begin{lem} \label{ambarbound}
With $L$  chosen as in \eqref{s1:L}, $A_m$ Hurwitz with constants $\sigma$ and $a$ as defined in \eqref{s1:ambound}, $P$ in \eqref{s1:lyap} satisfies
\begin{align}\label{s1:pbound}
 \text{(i)} & &\norm P \leq& \frac{m^2}{\sigma + 2\ell}  & &\\ 
 \text{(ii)} & &\min_i \lambda_i{(P)} \geq &\frac{1}{2(s+\ell)}& & \label{minP}
\end{align}
where $m=(1+4\varkappa)^{n-1} \text{ and } \varkappa \triangleq \frac a \sigma$.\end{lem}
\begin{proof} See Appendix \ref{pf:ambound}.\let\IEEEQED\relax\end{proof}

\subsection{Projection Algorithm}
Before we evaluate the benefits of closed--loop reference models, we introduce a modification in the adaptive law to ensure robustness properties. 
\begin{assm}\label{asm:thetas}
A  known $\theta^*_{max}$ exists such that $\norm{\theta^*}\leq \theta^*_{max}$.\end{assm}
The projection based adaptive law, which replaces \eqref{s1:law}, is given by
\be\label{s2:adaplaw}
\dot \theta(t) = \text{Proj}_\Gamma \left(\theta(t),-x e^T P b ,f \right)  
\ee
where the $\Gamma$--projection function, $\text{Proj}_\Gamma$, is defined as in Appendix \ref{ap:proj} and $f$ is a convex function given by
\be\label{eq:f}
f(\theta;\vartheta,\varepsilon) = \frac{\norm \theta^2 -\vartheta^2}{2 \varepsilon \vartheta - \varepsilon^2}
\ee
where $\vartheta$ and $\varepsilon$ are positive constants chosen as $\vartheta = \theta^*_\text{max}$ and $\varepsilon>0$.

\begin{defn}
Using the design parameters of the convex function $f(\theta;\vartheta,\varepsilon)$ we introduce the following definitions
\be\begin{split}\label{tmaxdef}
\theta_\text{max}&\triangleq  \vartheta+\varepsilon \text{ and }\\
\tilde\theta_\text{max} &\triangleq 2\vartheta+\varepsilon.
\end{split}\ee
\end{defn}

%
%
%


\subsection{Convergence of the Adaptive System}

\begin{thm}\label{thm:t2}
Let Assumptions \ref{asm:m} and \ref{asm:thetas} hold. Consider the adaptive system defined by the plant in \eqref{s1:sys} with the reference model in \eqref{s1:ref}, the controller in \eqref{s1:u}, the adaptive tuning law in \eqref{s2:adaplaw} and $L$ and $\Gamma$ as in \eqref{s1:L}-\eqref{s1:G}. For any initial condition in $e(0)\in \Re^n$, and $\theta(0)$ such that $\norm{\theta(0)}\leq\theta_\text{max}$, $e(t)$ and $\theta(t)$ are uniformly bounded for all $t\geq0$ and converge exponentially to a set $\mathcal E$ through $V$ in \eqref{can} as
\be\label{V_1}
\dot V \leq - \alpha_1 V + \alpha_2
\ee 
where
\be\label{alpha1}
 \alpha_1  \triangleq \frac{\sigma+2\ell}{m^2}\text{ and }  \alpha_2  \triangleq\frac {\sigma+2\ell} {m^2 \gamma}  \tilde\theta_\text{max}^2,
\ee
and\ben
\mathcal E \triangleq \left\{ (e,\tilde\theta) \left| \norm e^2 \leq \beta_1 \tilde\theta_\text{max}^2, \ 
  \norm{\tilde \theta} \leq \tilde\theta_\text{max }\right.\right \}\een
with 
\be\label{b1}
\beta_1 = 2 \frac{s+\ell}{ \gamma}.
\ee
\end{thm}
\begin{proof}See Appendix \ref{pf:thm2}. \let\IEEEQED\relax \end{proof}

%

\subsection{CRM free design parameters}
In this section it is argued that the free design parameters for the adaptive system are $\ell$ and the ratio of $\ell$ and $\gamma$. It is clear from \eqref{V_1} that the rate of exponential decay of the Lyapunov function is solely a function of the slowest eigenvalue of $\sigma$ and $\ell$. The term $\sigma$ is defined by the open--loop reference model jacobiam $A_m$ in \eqref{s1:ref} and cannot be independently increased without also increasing the bound on $\theta^*$, due to the fact that $A_m= A_p+b\theta^{*T}$ from \eqref{s1:match}.  Therfore, from \eqref{alpha1} it is clear that while increasing $\sigma$ increases $\alpha_1$ it also increases the size of of the compact set $\mathcal E$ that the model following error is exponentially converging to. The Luenberger gain $\ell$, in contrast, does not affect the matching condition for the adaptive system and therefore increasing $\ell$ does not result in a high--gain matching condition. It is seen from \eqref{b1} that while increasing $\ell$ may result in a larger $\alpha_1$ which is desirable, it also increases the $\alpha_2$ and $\beta_1$. At first glance this seems undesirable. However, it is important to note that $\ell$ is inversely proportional to the bounds for  $P$ in \eqref{s1:pbound} and \eqref{minP}. Upon further inspection of the tuning law in \eqref{s2:adaplaw}, when Projeciton is not active the adaptive law reduces to
\ben
\dot\theta= -\gamma  x e^T P b,
\een
and therefore, the only fixed parameters that control the rate of adaptation are $\gamma, P(\ell)$ and $b$. Given that $b$ is fixed it is ignored and thus the learning rate is a function of $\gamma \norm P $. Therefore, if $\gamma$ and $\ell$ are increased at the same rate, the effective rate of change in the adaptive tuning law will remain the same, which illustrates that the ratio of $\ell$ and $\gamma$ is an important design parameter. 
We therefore introduce $\rho$, an {\em effective learning rate}, as
\be\label{rhos}
\rho = \frac{\gamma}{\sigma + \ell}.
\ee

Using the definition of $\rho$ the bounds from Theorem \ref{thm:t2} can be rewritten as
\be\label{newbounds}
\begin{split}
\alpha_1  =& \frac{\sigma+2\ell}{m^2}\text{ and }  \alpha_2  \triangleq\frac 2 {m^2 \rho}  \tilde\theta_\text{max}^2, \\
\beta_1 =&  \frac{2s }{\sigma}\frac{1}{ \rho}.
\end{split}
\ee
This reparameterization in terms of $\rho$ and $\ell$ is used for discussing the transient performance.

\section{Transient Performance of CRM--adaptive systems}\label{section:perf}

In the following subsections we derive the transient properties of the CRM-adaptive systems. Five different subsections are presented, the first of which quantifies the Euclidean and the $\mathcal L_2$--norm of the tracking error $e$. In the second subsection we compute the same norms for the parameter derivative $\dot\theta(t)$. In both cases, we show that the $\mathcal L_2$--norms can be decreased by increasing $\ell$. In the third theorem, we address the performance of the true error $e^o$ and show its dependence on $\ell$. In the fourth subsection, we define our metric for transient performance in terms of a truncated $\mathcal L_2$ norm of the rate of control effort. The last subsection compares ORM and CRM adaptive systems using these metrics.
%

\subsection{Bound on $e(t)$}
\begin{thm} \label{ebound}Let Assumptions \ref{asm:m} and \ref{asm:thetas} hold. Consider the adaptive system defined by the plant in \eqref{s1:sys} with the reference model in \eqref{s1:ref}, the controller in \eqref{s1:u}, the adaptive tuning law in \eqref{s2:adaplaw} and $L$ and $\Gamma$ as in \eqref{s1:L} and \eqref{s1:G}.
\begin{align}
\label{e2e}\norm {e(t)}^2 &\leq  \kappa_1 \norm{e(0)}^2 \exp \left(-\frac{\sigma+2\ell}{m^2} t \right) + \frac{\kappa_2}{\rho} {\tilde\theta_\text{max}}^2\\
\label{el2}\norml{e(t)}{2}^2 &\leq \frac{1}{\sigma+ \ell} \left(m \norm{e(0)}^2 +  \frac {1} {\rho}{\norm{\tilde\theta(0)}}^2 \right)
\end{align} where $\kappa_i$, $i=1,2$ are independent of $\rho$ and $\ell$.\end{thm}
\begin{proof}
see Appendix \ref{eboundp}.\end{proof}



\subsection{Bound on $\dot\theta(t)$}

In addition to $\norm{e(t)}_{L_2}$ we explicitly compute upper bounds for $\norm{\dot\theta(t)}$ and $\norm{\dot{\theta}(t)}_{L_2}$. From the definition of $\dot\theta(t)$ in \eqref{s2:adaplaw}, it follows that
\ben
\norm{\dot\theta(t)} \leq \norm \Gamma \norm P  \norm b \norm{x(t)} \norm{e(t)}.
\een
We note that $x(t)= e(t) + x_m(t)$ and from \eqref{s1:ref} and \eqref{exp2} that
\be
\label{axm44}
\begin{split}
\norm{x_m(t)}\leq  &x_m(0) m  \exp\left( -\tfrac{\sigma}{2} t\right) \\ &+ m \int_0^t  \exp\left( -\tfrac{\sigma}{2} (t-\tau) \right) \left( \ell \norm{ e}+\norm b \norm r\right) d \tau
\end{split} \ee 
Using the bound for $\norm{e(t)}_{L_2}$ from \eqref{s2:adaplaw} and the Cauchy--Schwartz inequality, we simplify \eqref{axm44} as
\be\label{uhm}
\norm{x_m(t)} \leq  x_m(0) m  \exp\left( -\tfrac{\sigma}{2} t\right)+ \frac{\ell m}{\sqrt\sigma} \norml{e(t)}{2}+ \frac{ r_0 2 \norm{b} m} {\sigma}. 
\ee
The above bounds make the following theorem possible.

\begin{thm} \label{tbound}Let Assumptions \ref{asm:m} and \ref{asm:thetas} hold. Consider the adaptive system defined by the plant in \eqref{s1:sys} with the reference model in \eqref{s1:ref}, the controller in \eqref{s1:u}, the adaptive tuning law in \eqref{s2:adaplaw} and $L$ and $\Gamma$ as in \eqref{s1:L} and \eqref{s1:G}.
\begin{align}
\begin{split}\label{t1bound}
\norm{\dot\theta(t)}   \leq &\rho \exp{ \left(-\tfrac{\sigma+2\ell}{2m^2}t\right)} \left[ a_1 +\sqrt{\ell}\left( a_2 + a_3 \sqrt{\tfrac{1}{\rho}}\right)  \right] \\ &+ \sqrt{\rho} \exp{\left( -\tfrac{\sigma}{2} t\right)} a_4 +\sqrt{\tfrac 1 \rho} \exp{ \left(-\tfrac{\sigma+2\ell}{m^2}\right)} a_5 \\
&+ \sqrt{\ell\rho} a_6 + \sqrt{\ell} a_7 + \rho a_8
\end{split} \\ \label{tbound2}
\norml{\dot\theta(t)}{2} ^2 \leq  &  \rho^2 \nu(\rho) \left(\frac{b_1}{\sqrt{\sigma+\ell}} + \sqrt{\nu(\rho)} b_2 + \frac{b_3}{\sqrt{\sigma+\ell}} \right)^2
\end{align}
where $\nu(\rho)=m \norm{e(0)}^2 +  \frac {1} {\rho}{\norm{\tilde\theta(0)}}^2$, and the $a_i$ and $b_i$ are independent of $\rho$ and $\ell$.
\end{thm}
\begin{proof} see Appendix \ref{tboundp}.\end{proof}


\subsection{Bound on $e^o(t)$}
Here, we derive a bound on  the true error $e^o(t)$ defined in \eqref{eo}.

\begin{thm} \label{xmbound} Let the assumptions from Theorem \ref{tbound} hold. The difference between the open--loop reference model and the closed loop reference model satisfy the following bound
\be
\norm{e^o(t)} \leq \norm{e(t)} +  \sqrt{\frac \ell \sigma} m\sqrt{\nu(\rho)}.
\ee
\end{thm}
\begin{proof}see Appendix \ref{xmboundp}\end{proof}


\subsection{Bound on $\dot u (t)$}
We now derive a final transient measure of the CRM--adaptive system that pertains to $\dot u$. This is chosen as the transient performance metric because the rate of change of the control authority requested by the controller directly affects the robustness of the system to unmodelled dynamics and actuator rate limits. Before the bounds are derived, several variables must be defined.

\begin{defn}\label{allt}
Let time-constants $\tau_1(\ell)$, $\tau_2$ be defined as \be\label{taudef1}
\tau_1(\ell)= \frac{2m^2}{ \sigma  +2\ell} \text{ and }
\tau_2=\frac 2 \sigma
\ee
Let constants $a_\theta$ and $\delta_1(\ell,N)$ be defined as
\be\label{es}\begin{split}
a_\theta \triangleq& a+ \norm b \tilde\theta_\text{max},\\
\delta_1(\ell,N)=&\exp{(a_\theta N \tau_1(\ell))}-1.  \end{split}\ee 
where $N>0$, and three intervals of time\be
\begin{split}
\mathbb T_1 &= [0,N\tau_1)\\ 
\mathbb T_2 &= [N\tau_1,T_1) \\
\mathbb T_3 &=[T_1,\infty)
\end{split}
\ee
where $T_1\triangleq\max\{N\tau_2, T(\epsilon,-\ell I_{n\times n})\}$ and $T(\epsilon,-\ell I_{n\times n})$ is defined in Corollary \ref{Tdef}.\end{defn}  
\begin{rem}
$t_1(\ell)$ is a time constant associated with the exponential decay of $\norm{e(t)}$ which is derived from the upper bound on $V$ from \eqref{V_1} and $\tau_2$ is the time constant associated with $A_m$ in \eqref{s1:ref}. $a_\theta$ is a positive scalar that upper bounds the open--loop eigen values of $A_p$ from \eqref{s1:ref} and $\delta_1(\ell)$ will be used in the following Lemma to formally define our time scale separation condition. The time interval $\mathbb T_1$ is the time interval over which $\norm{e(t)}$ decays by $N$  time constants, $\mathbb T_3$ is the asymptotic time scale for $e(t)$ and $\mathbb T_2$ is an intermediate time interval. We note that $T_1$ exists but is unknown. 
\end{rem}

\begin{lem}\label{tscalesep}
For any $N>0$ an $\ell^*$ exists such that
\begin{enumerate}[(i)]
\item $\label{d} \delta_1(\ell^*,N)<\delta$  where $0<\delta\leq1$.
 \item $\tau_1(\ell^*)\leq \tau_2$.
\end{enumerate}
\end{lem}

\begin{rem} The condition Lemma \ref{tscalesep} (i) defines the time scale separation condition. Recall that $\tau_1$ is the time scale associated with $e(t)$ and $a_\theta$ is an upper bound on the uncertain open--loop eigen values of the plant. When $\ell\geq \ell^*$ we are able to show that at $t_N=N\tau_1$, $e(t_N)$ has exponentially decade by $N$ time constants, while $x(t_N)$ has not deviated far from $x(0)$.
\end{rem}

\begin{assm}\label{asst1}
$\exists r_0,r_1>0$ s.t. $ \abs{r(t)}\leq r_0$, ${\abs{\dot r(t)} \leq r_1}$.
\end{assm}
\begin{rem}
The bound on $\dot r(t)$ is needed so that $\dot u(t)$ is well defined. The analysis techniques that follow in proving Theorem \ref{thmudot} will still hold for reference inputs with discontinuities. The metric for transient performance however would change from $\dot u$ to $\frac{d}{dt} \left(\theta^T(t)x(t)\right)$.
\end{rem}

\begin{assm}\label{assxm0}
For ease of exposition we will assume that $ x_m(0)=0.$ \end{assm} We note that the same analysis holds for $x_m(0)$ with addition of exponentially decaying terms proportional to $x_m(0)$.


\begin{thm}\label{thmudot}
Let Assumptions 1--4 hold. Given arbitrary initial conditions in ${x(0)\in\Re^n}$ and ${\norm{\theta(0)}\leq \theta_\text{max}}$, for any  ${\epsilon>0}$, ${N>0}$ and ${\ell\geq\ell^*}$ , $\dot u$ satisfies the following inequalities:
\be\label{ub1}
\begin{split}
 \sup_{t\in\mathbb T_i}\abs{\dot u(t)} \leq &\frac{m^2\gamma}{\sigma+2\ell}\norm b G_{e,i} G_{x,i}^2  \\ &+ \theta_\text{max}\left(a_\theta G_{x,i}+r_0\right) + r_1
\end{split}
\ee for $i=1,2,3$, where
\be\label{GS}\begin{split}
G_{x,1}\triangleq & (1+\delta_1) \norm{e(0)} +\frac{\delta_1 \norm b }{a_\theta}r_0\\
G_{e,1}\triangleq & \sqrt{\kappa_1}\norm{e(0)}+\sqrt{\frac{\kappa_2}{\rho}}\tilde\theta_\text{max} \\
G_{x,2}\triangleq &\kappa_3 \norm{e(0)} + \left(1+\kappa_4 \ell\right)  \sqrt{\frac{\kappa_2}{\rho}}  \tilde\theta_\text{max} +\kappa_5 r_0\\
G_{e,2}\triangleq &\sqrt{\kappa_1}\norm{e(0)}\epsilon_1+\sqrt{\frac{\kappa_2}{\rho}}\tilde\theta_\text{max}\\
G_{x,3}\triangleq &\kappa_6\norm{e(0)}+\epsilon+\left(1+\kappa_4 \ell\right)  \sqrt{\frac{\kappa_2}{\rho}}  \tilde\theta_\text{max} +\kappa_5 r_0\\
G_{e,3}\triangleq & \epsilon
\end{split}
\ee
where ${\epsilon_1\triangleq\exp(-N)}$ and the $\kappa_i$ are independent of $\rho$ and $\ell$, and $N\geq 3$
\end{thm}
\begin{proof}see Appendix \ref{tan}.  \let\IEEEQED\relax \end{proof}

\begin{rem}
There are two ``small'' terms in the above analysis. $\epsilon_1$ is determined by the number of time constants $N$ of interest. $\epsilon$ is free to choose and from Corollary \ref{Tdef} proves the existence of a finite $T$ and is used to define when  $\mathbb T_3$ begins.  \end{rem}

From Theorem \ref{thmudot}, it follows that
\be\label{uox}\begin{split}
\sup_{t\in\mathbb T_1}\abs{\dot u(t)} \leq & c_1\rho + c_2 \sqrt\rho + r_1\\ 
  \sup_{t\in\mathbb T_2}\abs{\dot{u}(t)} \leq & \sqrt \rho c_3 +  (1+c_{4} l) c_5  +  \sqrt{\frac 1 \rho}  (1+c_{4} \ell)^2   c_6 \\
   &+ \epsilon_1 \mathfrak{L}_1 (\rho,\ell,\sqrt \rho, \ell\sqrt\rho, \ell^2) + r_1\\
    \sup_{t\in\mathbb T_3}\abs{\dot{u}(t)} \leq &  \sqrt{\frac 1 \rho} (1+c_4\ell )c_7 +c_8\\ &+\epsilon \mathfrak{L}_2 (\rho,\ell,\sqrt \rho, \ell\sqrt\rho, \ell^2,\epsilon_1) + r_1
 \end{split}
\ee
where $c_i>0$, $i=1$ to $8$  are independent of $\ell$ and $\rho$, $\mathfrak L_1(\cdot)$ and $\mathfrak L_2(\cdot)$ are globally lipschitz with respect to their arguments.
The inequalities in \eqref{uox} lead us to the following three main observations (see Figure \ref{fig:udot})
\begin{enumerate}[(\bf {A}1)]
\item Over $\mathbb T_1$, $\abs{\dot u(t)}$ is bounded by a linear function of $\rho$ and $\sqrt{\rho}$, 
\item Over $\mathbb T_2$, $\abs{\dot u(t)}$ is bounded by a linear function of $\sqrt \rho,\ell,\sqrt {\frac{1}{\rho}},\ell \sqrt {\frac{1}{\rho}}$ and $\ell^2\sqrt {\frac{1}{\rho}} $
\item Over $\mathbb T_3$, $\abs{\dot u(t)}$ is bounded by a linear function of $\sqrt {\frac{1}{\rho}}$ and $\ell \sqrt {\frac{1}{\rho}}$
\item $\tau_1$ decreases with $\ell$.
\end{enumerate}

\begin{rem} The main idea used for the derivation of the bounds in Theorem \ref{thmudot} is time--scale separation of the error dynamics decay, and the worst case open--loop eigenvalues of the uncertain plant. The most important point to note is that $\tau_1$ can be made small by choosing a large $\ell$. There is a penalty, however, in choosing a large $\ell$, as the bound $G_{x,2}$ increases linearly with $\ell$. Therefore, after choosing an $\ell$ which satisfies the time scale separation as needed in Lemma \ref{tscalesep}, a $\rho$ (which through \eqref{rhos} defines a choice for $\gamma$) can be chosen such that the integral in the following theorem is minimized. \end{rem}

\begin{thm}\label{thmpopt}
There exist optimal $\rho$ and $\ell$ such that
\be\label{optp}
\left(\rho_\text{opt},\ell_\text{opt}\right)=\argmin_{{\rho>0} \atop{\ell\geq\ell^*}} \norm{\dot u(\rho,\ell)}_{L_2, \tau}
\ee
for any $0 < \tau < T_1$.
\end{thm}
\begin{proof}
$\norm{\dot u(\rho,\ell)}_{L_2, \tau}$ is continuous with respect to $\rho$ and $\ell$ where $\rho$ and $\ell$ appear in the numerator of  \eqref{uox} and are positive. Therefore, $\rho_\text{opt}$ and $\ell_\text{opt}$ exist and are finite.
\end{proof}

$\tau$ in Theorem \ref{thmpopt} denotes the interval of interest in the adaptive system where the transient response is to be contained. Given that $T$, and therefore $T_1$ is a function of $\ell$, \eqref{optp} can only be minimized over $\mathbb T_1 \cup \mathbb T_2$. From the authors definition of smooth transient performance in the beginning of this section choosing $\rho_{opt}$ and $\ell_{opt}$ will guarantee smooth transient performance.

\subsection{Comparison of CRM and ORM-adaptive systems}\label{sqwe}
The bounds on  $e(t)$ and the $\mathcal L_2$--norm of $\dot\theta$ directly show that CRM--adaptive systems lead to smaller $e(t)$ than with the ORM which are obtained by setting $\ell=0$ in \eqref{e2e} and \eqref{el2}. However, the same cannot be said for either $e^o$ or for the Euclidean norm of $\dot\theta$; for a non-zero $\ell$, the bound on $e^o$ is larger than that of $e$.  This indicates that there is a trade-off between fast transients and true tracking error. The signal that succinctly captures this trade off is $\dot u$, whose behavior is captured in detail using the time intervals $\mathbb T1$, $\mathbb T2$, and $\mathbb T3$. We also showed in Theorem \ref{thmpopt} that this trade-off can be optimized via a suitable choice of $\ell$ and $\rho$. In what follows, we compare this optimized CRM with ORM and show that the former is clearly better than the latter.

\begin{defn}The following two time constants 
\be
\tau_2^\prime\triangleq  \tau_1(0)=\frac{2m^2}{\sigma} \text{ and } \tau_1^* = \tau_1(\ell^*) 
\ee
are used to describe the three time intervals that will be used in the analysis of $\dot u$ for the ORM case
\be
\begin{split}
\mathbb T_1^\prime &= [0,N\tau_1^*)\\ 
\mathbb T_2^\prime &= [N\tau_2^\prime,T_1^\prime)\\
\mathbb T_3^\prime &= [T_1^\prime,\infty). \end{split}
\ee
where $T_1^\prime\triangleq\max\{N\tau^\prime_2, T(\epsilon,0)\}$ where $T(\epsilon,0)$ is from Corollary \ref{Tdef}.\end{defn}

\begin{figure}[t!]
\centering
\psfrag{u}[cc][cc][.9][90]{$\abs{\dot u}$}
\psfrag{t1}[cc][cc][.7]{$N\tau_1^*$}
\psfrag{t2}[cc][cc][.7]{$N\tau_2$}
\psfrag{t3}[cc][cc][.7]{$N\tau_1$}
\psfrag{t5}[cc][cc][.9]{$\mathbb T_1$}
\psfrag{t6}[cc][cc][.9]{$\mathbb T_2$}
\psfrag{b1}[cc][cc][.7]{$T_1^\prime$}
\psfrag{b2}[cc][cc][.7]{$T_1$}

\psfrag{a1}[cc][cc][.9]{$\mathbb T_3^\prime$}

\psfrag{a2}[cc][cc][.9]{$\mathbb T_3$}

\psfrag{t7}[cc][cc][.9]{$\mathbb T_3$}
\psfrag{t8}[cc][cc][.9]{$\mathbb T_1^\prime$}
\psfrag{t9}[cc][cc][.9]{$\mathbb T_2^\prime$}

\psfrag{t4}[cc][cc][.7]{$N\tau_2^\prime$}
\psfrag{t}[cc][cc][.9]{$t$}
\psfrag{d}[cc][cc][.9]{$\Delta$}
\psfrag{e1}[cl][cl][.9]{$\ell=\ell^*$}
\psfrag{e2}[cl][cl][.9]{$\ell=0$}
\psfrag{g}[cl][cl][.9]{$O(g_l e(0))$}
\psfrag{s}[cl][cc][.9]{$O\left(e(0)\frac{\sqrt{q\gamma}}{\sigma+2\ell}\right)$}
\includegraphics[width=2.8in]{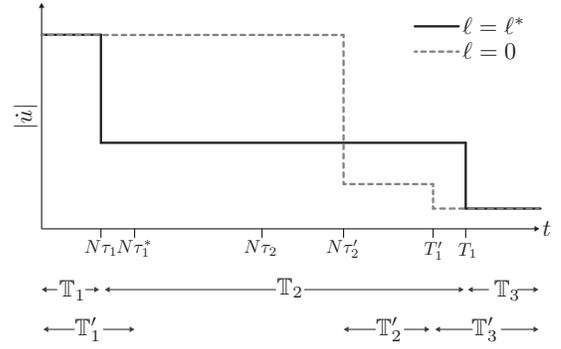}
\caption{Transient bounds for $\dot u$.}\label{fig:udot}
\end{figure}

As in Definition \ref{allt}, here too, $T$ exists but is unknown. While these periods for both CRM and ORM are indicated in Figure \ref{fig:udot}, one cannot apriori conclude if $T_1$ is greater than or smaller than $T_1^\prime$. The time instants indicated as in Figure \ref{fig:udot} are meant to be merely sketches.

\begin{prop}\label{propl0}
Let
\be
\rho_0 \triangleq  \frac{\gamma}{\sigma}.
\ee
For the adaptive system with the classical MRAC given by Eqs \eqref{s1:sys},  \eqref{s1:ref}, \eqref{s1:u}, \eqref{s2:adaplaw}--\eqref{eq:f}  and \eqref{s1:L}--\eqref{s1:G} with $\ell=0$, it can be shown that
\be\label{uox2}\begin{split}
\sup_{t\in\mathbb T_1^\prime}\abs{\dot u(t)} \leq & \rho_0 d_1 + \sqrt {\rho_0} d_2  + r_1,     \\
 \sup_{t\in\mathbb T_2^\prime}\abs{\dot{u}(t)} \leq & \sqrt {\rho_0} d_3  +d_4   + \sqrt{\frac 1 {\rho_0}}  d_5 + \epsilon_1 \mathfrak{M}_1 (\rho_0,\sqrt{\rho_0}) + r_1  \\
    \sup_{t\in\mathbb T_3^\prime}\abs{\dot{u}(t)} \leq &  \sqrt{\frac 1 {\rho_o}} d_6 +d_7+\epsilon \mathfrak{M}_2 (\rho_o,\sqrt{\rho_0}) + r_1
 \end{split}
   \ee
$d_i>0$, $i=1$ to $7$  are independent of $\rho_0$, and $\mathfrak M_1(\cdot)$ and $\mathfrak M_2(\cdot)$ are globally lipschitz with respect to their arguments
\end{prop}
The proof of Proposition \ref{propl0} follows the same steps as in the proof of Theorem \ref{thmudot} and is therefore omitted.

The bounds in \eqref{uox2} indicate that in the classical ORM, one can only derive a bound for $\dot u$ over the period $\mathbb T_1^\prime$, $\mathbb T_2^\prime$ and $\mathbb T_3^\prime$. Unlike the CRM case, the procedure in Appendix \ref{tan} cannot be used to derive satisfactory bounds for $\dot u$ over $[N\tau_1^*,N\tau_2^\prime)$. It also can be seen that unlike the CRM case, $\tau_2^\prime$ is fixed and cannot be changed with $\ell$. These points are summarized below. \begin{enumerate}[(\bf{B}1)]
\item Over $\mathbb T_1^\prime$, $\abs{\dot u(t)}$ is bounded by a linear function of $\rho_0$ and $\sqrt{\rho_0}$
\item Over $\mathbb T_2^\prime$, $\abs{\dot u (t)}$ is bounded by a linear function of $\sqrt {\rho_0}$ and $\sqrt {\frac{1}{\rho_0}}$
\item Over $\mathbb T_3^\prime$, $\abs{\dot u (t)}$ is bounded by a linear function of $\sqrt {\frac{1}{\rho_0}}$
\item $\tau_2^\prime$ is fixed and unlike $\tau_1$, can not be adjusted.
\end{enumerate}

We now compare the bounds on $\dot u$ using observations (A1)--(A3) and (B1)--(B3). In order to have the same basis for comparison, we assume that $\gamma$, $\sigma$, and $\ell$ are such that $\rho=\rho_0$ and that both CRM-- and ORM--adaptive systems start with the same bound at $t=0$. As noted above, a tight bound cannot be derived for the ORM-based adaptive system over $[N\tau_1^*,N\tau_2^\prime)$.  In the best scenario, one can assume that this bound is no larger than that over $[0,N\tau_1^*]$. This allows us to derive the bounds shown in Figure \ref{fig:udot}. The main observations that one can make from this figure are summarized below:
\begin{itemize}
\item  Even though at time $t = 0$, both the ORM and CRM have the same bound, since $\tau_1$ can be made much smaller than $\tau_2^\prime$, this bound is valid for a much shorter time with the CRM-system than in the ORM--system. This helps us conclude that the initial transients can be made to subside much faster in the former case than the latter, by suitably choosing $\ell$. 
\item The bound on $\dot u$ for $\mathbb T_2$ with the CRM--adaptive system is however linear in powers of $\ell$ and hence can be larger than the bound on $\dot u$ with the ORM-adaptive system over $\mathbb T_2^\prime$. 
\item The above observations clearly illustrate, if the cost function $U(N\tau_2^\prime;\rho,\ell)$ is minimized then the CRM system will have smoother transients than the ORM. Then, at larger times the error dynamics will asymptotically converges to zero.
\end{itemize}

\begin{figure}[t!]
\centering
\psfrag{u}[cc][cc][.9]{$u$}
\psfrag{u1}[cc][cc][.9]{$\abs{\dot u(t)}$}
\psfrag{t}[cc][cc][.8]{$t$}
\psfrag{r}[cc][cc][.8]{$x_m$}
\psfrag{openloop}[cl][cl][.8]{open--loop}
\psfrag{closedloop}[cl][cl][.8]{closed--loop}
\psfrag{extrascenario1}[cl][cl][.8]{$\ell$=0 \, $\rho=100$}
\psfrag{extrascenario2}[cl][cl][.8]{$\ell$=10 $\rho=100$}
\psfrag{extrascenario3}[cl][cl][.8]{$\ell$=10 $\rho=1$}
\includegraphics[width=2.8in]{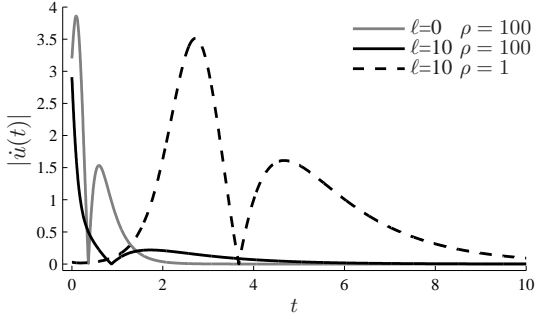}
\caption{Plot of $\abs{\dot u(t)}$.}\label{fig:ta1}
\end{figure}

\subsection{Water--Bed Effect}

The discussions in the preceding sections clearly show that CRM-adaptive systems introduce a trade--off: a fast convergence in $e(t)$ with a reduced $\norml{\dot\theta(t)}{2}$ occurs at the expense of an increased $e^o(t)$. While an optimal choice of $\rho$ and $\ell$ can minimize this trade--off, it also implies that a badly chosen $\ell$ and $\rho$ can significantly worsen the adaptive system performance in terms of $e^o(t)$ and $\dot u(t)$. We denote this as the water--bed effect and illustrate it through a simulation. Consider a first-order plant with a single unknown parameter, whose values are identical to the example in Section \ref{simstud1} over the first ten seconds. Figure \ref{fig:ta1} shows the behavior of $\dot u(t)$ for the ORM, the optimized CRM, and a poorly chosen CRM. The plots clearly show the water--bed effect for the last case and the improved performance of the optimized CRM over the ORM. The free design parameters are also shown in the figure.

\section{Robustness of CRM to Time--Varying Uncertainties and Disturbances}
We now evaluate the CRM--adaptive system in the presence of perturbations due to time-varying parameters and disturbances. Consider the uncertain Linear Time Varying system
\be\label{s2:sys}
\dot x = A_p(t) x(t) + b u  + d(t)
\ee
where $d(t)$ is a bounded disturbance and $A_p(t)$ is time varying with a bounded time--derivative. It is assumed that a time-varying vector $\theta^*(t)$ exists such that 
\be
A_m = A_p (t)+ b \theta^{*T}(t),
\ee
and $\theta_d^*$, $\theta_{max}$ exist where $\norm{\dot \theta^*(t)}\leq{\theta_d^*}$ and $\norm{\theta^*(t)}\leq \theta_{max}.$

\begin{thm}\label{lem:s2}
With Assumptions \ref{asm:m} and \ref{asm:thetas}, consider the adaptive system defined by the plant in \eqref{s2:sys} with the reference model in \eqref{s1:ref}, the controller in \eqref{s1:u}, the adaptive tuning law in \eqref{s2:adaplaw} and $L$ and $\Gamma$ as in \eqref{s1:L}-\eqref{s1:G}. For any initial condition in $e(0)\in \Re^n$, and $\theta(0)$ such that $\norm{\theta(0)}\leq\theta_\text{max}$, $e(t)$ and $\theta(t)$ are uniformly bounded for all $t\geq0$ and the Lyapunov candidate in \eqref{can} converges exponentially to a set $\mathcal E$ as\be
\dot V \leq - \alpha_3 V + \alpha_4
\ee 
where $ \alpha_3  \triangleq \frac{\alpha_1}{2}$,
\be\label{alpha4}
\alpha_4  \triangleq  \frac {\sigma+2\ell} {2m^2 \gamma} \tilde\theta_\text{max}^2 + \frac {2} \gamma \theta^*_d\tilde \theta_\text{max}  +2 \left(\frac{m^2}{\sigma+2\ell}\right) ^2 \norm {d(t)}^2,
\ee
and\ben
\mathcal E \triangleq \left\{ (e,\tilde\theta) \left| \norm e^2 \leq \beta_1  \tilde\theta_\text{max}^2 + \beta_2\theta^*_d\tilde\theta_\text{max} +\beta_3 \norm d^2,\right.
  \norm {\tilde \theta} \leq \tilde\theta_\text{max} \right\}
\een
where 
\be\label{b2}
\beta_2 \triangleq \frac{8s m^2}{ \sigma \gamma }  \text{ and } \beta_3 \triangleq  \frac{4 s m^6}{\sigma (\sigma+\ell)^2}.
\ee
\end{thm}
\begin{proof}See Appendix \ref{thm34}. \let\IEEEQED\relax \end{proof}

\begin{rem}
From the above Theorem it is is shown that in the presence of disturbance and time--varying uncertainty, the rate of exponential decay of the lyapunov function is $\alpha_1/2$ where $\alpha_1$ was the rate of exponential decay for the system with constant uncertainty and no disturbances. The term that controls how the compact set $\mathcal E$ grows with the upper bound on the rate of change of the time-varying matched uncertainty is $\beta_2$ which is inversely proportional to $\gamma$. $\beta_3$ scales the size of the compact set in regard to the disturbance term $d$. Therefore, increasing $\gamma$ and $\ell$ in proportion so that  $\rho$ is constant constant decreases the size of the compact set $\mathcal E$.
\end{rem}

\subsection{Simulation Study} \label{simstud1}
For this study a scalar time varying system of the form in \eqref{s2:sys} is controlled where
\ben
A_p(t) = \begin{cases}   1 & 0 \leq t < 20 \\
                                   1+\frac{1}{4} (t-20) & 20 \leq t <24 \\
                                   2 & t \geq 24  \end{cases}, 
\een
$b=1$, and $d(t)$ is a deterministic signal used to represent a disturbance. Over the first 20 seconds $d(t)=0$. After 20 seconds $d(t)$ is generated from a Gausian distribution centered at 0 with a variance of 1, covariance of 0.1, deterministically sampled at 10 Hz with a fixed seed, and then passed through a saturation function with upper and lowers bounds of 0.2 and -0.2 respectively. The reference model to be followed is defined as
\be\label{ex:m4}
\dot x_m = - x_m + r + \ell(x-x_m)
\ee
with control input from \eqref{s1:u} and the update law for the adaptive parameter defined \eqref{s1:law},
where $\ell$ and $\rho$ are chosen as in Table \ref{t:scen}.

\begin{figure}[t]
\centering
\psfrag{e}[cc][cc][.9]{$e$}
\psfrag{x}[cc][cc][.9]{$x$}
\psfrag{t}[cc][cc][.8]{$t$}
\psfrag{r}[cc][cc][.9]{$x_m$}
\psfrag{openloop}[cl][cl][.8]{open--loop}
\psfrag{closedloop}[cl][cl][.8]{closed--loop}
\psfrag{closedloope0}[cl][cl][.8]{$e^o$: closed--loop}
\psfrag{c1}[cc][cl][.8]{Region 1}

\psfrag{c2}[cc][cl][.8]{Region 2}

\psfrag{c3}[cc][cl][.8]{Region 3}
\includegraphics[width=3.4in]{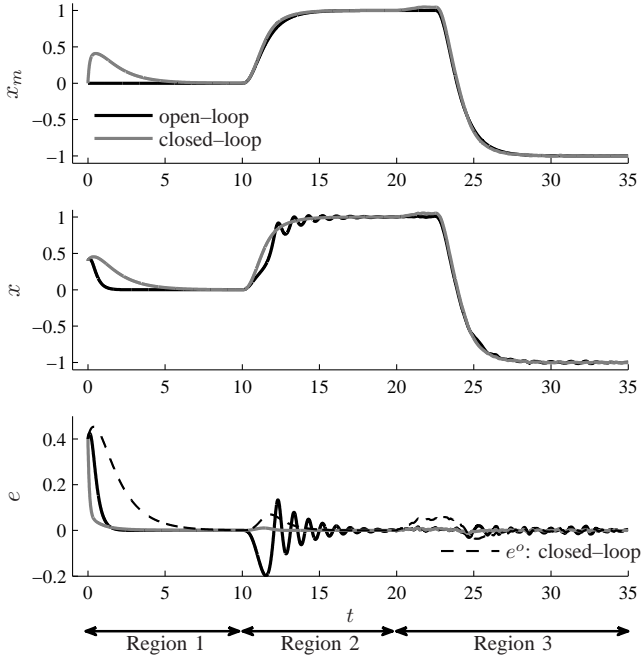}
\caption{(top) reference model trajectories $x_m$, (middle) state $x$, and (bottom) model following $e$.}\label{fig:e}
\end{figure}

\begin{figure}[h!]
\centering
\psfrag{u}[cc][cc][.9]{$u$}
\psfrag{u1}[cc][cc][.9]{$\dot u$}
\psfrag{t}[cc][cc][.8]{$t$}
\psfrag{r}[cc][cc][.8]{$x_m$}
\psfrag{h}[cc][cc][.9]{$\theta$}
\psfrag{thetastar}[cl][cl][.8]{$\theta^*(t)$}
\psfrag{openloop}[cl][cl][.8]{open--loop}
\psfrag{c1}[cc][cl][.8]{Region 1}

\psfrag{c2}[cc][cl][.8]{Region 2}

\psfrag{c3}[cc][cl][.8]{Region 3}
\psfrag{closedloop}[cl][cl][.8]{closed--loop}
\includegraphics[width=3.4in]{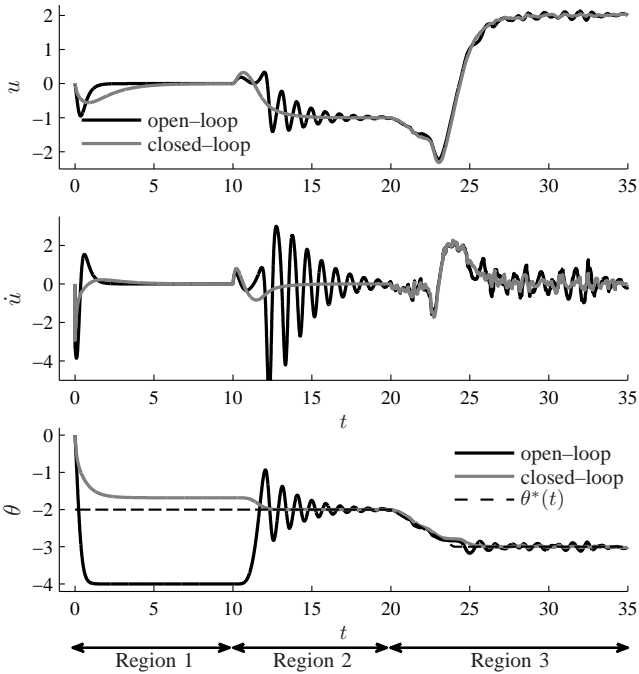}
\caption{(top) Control input $u$, (middle) rate of control input $\dot u$, and (bottom) adaptive parameter $\theta(t)$.}
\label{fig:u}
\end{figure}

\begin{table}[h!]
\caption{Test case free design parameters} \label{t:scen}
\centering
\begin{tabular}{ c c c}
Parameter & Open--Loop & Closed--Loop  \\ \hline
 $\ell$ & 0 & 10 \\
   $\rho$  & 100 & 100 \\ \hline
\end{tabular}  
\end{table}

The simulations have three distinct regions of interest, with Region 1 denoting the first 10 seconds, Region 2 denoting the 10 sec to 20 sec range, and Region 3 denoting the 20 sec to 35 sec range. In Region 1, the adaptive system is subjected to non--zero initial conditions in the state and the reference input is zero. At $t=10$ sec, the beginning of Region 2, a filtered step input is introduced. At $t=20$ sec In Region 3, time-variations in the plant parameter as well, disturbances and a filtered step input are introduced. Figures 2 and 3 illustrate the response of the CRM--adaptive system over 0 to 35 seconds, with $x_m$, $x$, and $e$ indicated in Figure 2, and $u$, $\dot u$, and $\theta$ indicated in Figure 3. In both cases, the resulting performance is compared with the classical adaptive system. The first point that should be noted is a satisfactory behavior in the steady-state of the CRM--adaptive controller. In particular, as can be seen from the latter half of Region 1, both $e$ and $e^o$ tend to zero as $t$ approaches 10 seconds. The same satisfactory trends are observed in regions 2 and 3 as well, underscoring the robustness property of CRM--adaptive control, which validates Theorem 2.
 
We also note yet another significant difference between the responses of CRM--adaptive control and the classical one, which pertains to the rate of control input $\dot u$. An examination of Regions 2 and 3 clearly illustrates that the control input is smoother for CRM--adaptive control.

\section{CMRAC}

We now return to CMRAC introduced in \cite{duarte1989combined} and \cite{slotine1989composite}. We will show that the introduction of a CRM in these adaptive systems not only ensures stability, but also enables stability with observer--based rather than state--based feedback. In addition, the use of a CRM in CMRAC enables the derivation of transient properties which could not be accomplished hitherto. Section A addresses stability of the CMRAC with CRM, denoted as CMRAC--C, in the scalar case. Section B extends the results from section A to higher order plants with states accessible and addresses transient properties of the CMRAC-C adaptive system. Section C introduces an additional feature of observer feedback. Denoting the underlying adaptive system as CMRAC--CO, it is shown that the resulting adaptive system has guaranteed stability properties and results in reduced error bounds in the presence of measurement disturbances with a zero mean property, which is corroborated through simulations in Section D. Section E contains extensions to higher order plants whose states are accessible.

\subsection{Stability of CMRAC--C}
We assume that the plant and reference model dynamics are given by Equations \eqref{s1:sys} and \eqref{s1:ref} with $A_m$ and $L=L_m$ satisfying Equations \eqref{s1:match} and \eqref{s1:ambar}. The control input is chosen as in \eqref{s1:u}
and the identifier dynamics are given by
\be\label{s8:o}
\dot x_i(t) = L_i (x_i(t) - x(t)) + (A_m- b\hat \theta^T(t) )x(t) + b u(t)
\ee
where $L_i$ is Hurwitz. The error dynamics are now given by
\be
\begin{split}\label{s8:e}
\dot e_m(t) =& (A_m+L_m) e_m +b \tilde \theta ^T (t)x\ \\
\dot e_i(t) =& L_i e_i -b \bar \theta ^T (t)x, \quad e_i=x_i-x
\end{split}
\ee
where ${\bar\theta(t)=\hat\theta(t)-\theta^*}$. For ease of exposition we choose
\be\label{s8:l}
L_m=L= -\ell  I_{n\times n} \text{ and } L_i=-(\sigma+\ell)I_{n\times n}.
\ee
The update laws for the adaptive parameters are then defined with the update law
\be\label{s8:law}\begin{split}
\dot \theta = \text{Proj}_\Gamma  (\theta(t),-x e_{m}^T P_m b,f)  - \eta I_{n\times n} \epsilon_\theta\\
\dot {\hat \theta} = \text{Proj}_\Gamma (\hat\theta(t), x e_i^T P_i b,f)+\eta I_{n\times n} \epsilon_\theta
\end{split}\ee
where $\epsilon_\theta = \theta(t) - \hat\theta(t)$, with $\Gamma$ chosen as in \eqref{s1:G}, ${\eta>0}$, and ${P_m=P}$ from \eqref{s1:lyap} and ${P_i=\frac{1}{2(\sigma+\ell)}I_{n\times n}}$.

\begin{thm}\label{cmracsat} Let Assumptions \ref{asm:m} and \ref{asm:thetas} hold. Consider the overall CMRAC--C specified by \eqref{s1:sys}, \eqref{s1:ref}, \eqref{s1:u},  \eqref{s8:o}, \eqref{s8:e} and \eqref{s8:law}. For any initial condition $e_m(0),e_i(0)\in \Re^n$, and $\theta(0)$ and $\hat\theta(0)$ such that $\norm{\theta(0)}\leq\theta_\text{max}$ and $\norm{\hat\theta(0)}\leq\theta_\text{max}$, it can be shown that $e_m(t)$, $e_o(t)$, $\theta(t)$ and $\hat\theta(t)$ are uniformly bounded for all $t\geq0$ and the function 
\be\label{s8:v}
V=   e_m^TP_m e_m + e_i^T P_i e_i + \tilde \theta^T\Gamma^{-1}\tilde\theta + \bar\theta^T \Gamma^{-1} \bar\theta
\ee 
converges exponentially to a set $\mathcal E$ as
\be\label{tbone3}
\dot V \leq -\alpha_1 V + 2 \alpha_2
\ee
where
\ben\begin{split}
\mathcal E \triangleq \left\{ (e_m,e_i,\tilde\theta,\bar\theta) \right | & \norm{e_m}^2  \leq \beta_4 \tilde\theta_\text{max}^2,  \norm{e_i}^2 \leq \beta_5  \tilde\theta_\text{max}^2 \\  & \left.  \norm {\tilde \theta} \leq \tilde\theta_\text{max},\ \norm {\bar \theta} \leq \tilde\theta_\text{max}\right\}
\end{split}\een
with 
\be\label{s4:b}
\beta_4 \triangleq \frac{4(s+l)}{\gamma} \text{ and } \beta_5\triangleq \frac{4(\sigma+\ell)}{\gamma}.
\ee
\end{thm}
\begin{proof}
see Appendix \ref{proofcmracsat}. \let\IEEEQED\relax  
\end{proof}

\begin{rem}
There is no appreciable difference between the CMRAC--C and CRM adaptive controller presented in Section \ref{sec:l}  in terms of stability and the bounds for the set $\mathcal E$.
\end{rem}

\section{Transient Properties of CMRAC-C}\label{section:perf2}

In the following subsections we derive the transient properties of the CMRAC--C adaptive system, similar to what was done in Section \ref{section:perf}. Two different subsections are presented, the first of which quantifies the Euclidean and the $\mathcal L_2$--norm of the tracking error $e$ and the second subsection, were the truncated $\mathcal L_2$ norm of the rate of control effort is presented. 

\subsection{Bound on $e_m(t)$ and $e_i(t)$}
\begin{thm} \label{eperf}Let Assumptions \ref{asm:m} and \ref{asm:thetas} hold. Consider the overall CMRAC--C specified by \eqref{s1:sys}, \eqref{s1:ref}, \eqref{s1:u},  \eqref{s8:o}, \eqref{s8:e} and \eqref{s8:law}. For any initial condition $e_m(0),e_i(0)\in \Re^n$, and $\theta(0)$ and $\hat\theta(0)$ such that $\norm{\theta(0)}\leq\theta_\text{max}$ and $\norm{\hat\theta(0)}\leq\theta_\text{max}$.
\begin{align}\begin{split}
\label{e2me}\norm {e_m(t)}^2 \leq& \kappa_7\left( \norm{e_m(0)}^2 +\norm{e_i(0)}^2\right) \exp \left(- \alpha_1 t \right) \\ &+ \frac{\kappa_8}{\rho} {\tilde\theta_\text{max}}^2\end{split} \\
\norm {e_i(t)}^2 \leq&\norm {e_m(t)}^2 \label{e2i} \\
\begin{split}
\label{eml2}\norml{e_m(t)}{2}^2 \leq& \frac{1}{\sigma+ \ell} \left(m^2 \norm{e_m(0)}^2 +\norm{e_i(0)}^2\right) \\& +\frac{1}{\sigma+ \ell} \left( \frac {1} {\rho} {\norm{\tilde\theta(0)}}^2+\frac {1} {\rho} {\norm{\bar\theta(0)}}^2 \right)\end{split}\\
\norml{e_i(t)}{2}^2 \leq& \norml{e_m(t)}{2}^2\label{eei2}
\end{align} where $\kappa_i$, $i=7,8$ are independent of $\rho$ and $\ell$.\end{thm}
\begin{proof}
see Appendix \ref{proofeperf}.\end{proof}

\subsection{Bound on $\dot u(t)$}

\begin{defn} The following three time intervals are used when exploring the transients of CMRAC--C
\be
\begin{split}
\mathbb T_1^{\prime\prime} &= [0,N\tau_1)\\ 
\mathbb T_2^{\prime\prime} &= [N\tau_1,T_1^{\prime\prime})\\
\mathbb T_3^{\prime\prime} &= [T_1^{\prime\prime},\infty)\\
\end{split}
\ee
where $T_1^{\prime\prime}=\max\{N\tau_2,T(\epsilon,-\ell I_{n\times n})\}$, with $T(\epsilon,-\ell I_{n\times n})$ following from the application of Barbalat Lemma to the adaptive system defined in Thereom \ref{cmracsat} for any $\epsilon>0$ (identical to Corollary \ref{Tdef}).\end{defn}

\begin{thm}\label{thmudot3}
Let Assumptions 1--4 hold. Given arbitrary initial conditions in $x(0)\in\Re^n$ and $\norm{\theta(0)}\leq \theta_\text{max}$, if $\ell\geq\ell^\prime$ the derivative $\dot u$ satisfies the following two inequalities:
\be
\begin{split}
 \sup_{t\in T_i^{\prime\prime}}\abs{\dot u(t)} \leq &\left(\frac{m^2\gamma}{\sigma+2\ell}\norm b G^{\prime\prime}_{e,i}G^{\prime\prime}_{x,i} +8\eta  \theta^2_\text{max}\right)G^{\prime\prime}_{x,i}   \\ & + \theta_\text{max}\left(a_\theta G^{\prime\prime}_{x,i}+r_0\right) + r_1\label{ub13}
\end{split}
\ee where
\be\label{GS3}\begin{split}
G^{\prime\prime}_{x,1}\triangleq & (1+\delta_1) \norm{e(0)} +\frac{\delta_1 \norm b }{a_\theta}r_0\\
G^{\prime\prime}_{e,1}\triangleq & \sqrt{\kappa_7}\left( \norm{e_m(0)} +\norm{e_i(0)}\right) + \sqrt{\frac{\kappa_8}{\rho}} {\tilde\theta_\text{max}} \\
G^{\prime\prime}_{x,2}\triangleq &\kappa_9 \left( \norm{e_m(0)} +\norm{e_i(0)}\right) +
\left(2+\kappa_{10} \ell \right) \sqrt{\frac{\kappa_8}{\rho}} {\tilde\theta_\text{max}} 
\\ &+\kappa_{11} r_0\\
G^{\prime\prime}_{e,2}\triangleq  & \sqrt{\kappa_7}\left( \norm{e_m(0)} +\norm{e_i(0)}\right) \epsilon_1 + \sqrt{\frac{\kappa_8}{\rho}} {\tilde\theta_\text{max}} \\
G^{\prime\prime}_{x,3}\triangleq &\kappa_{12} \left( \norm{e_m(0)} +\norm{e_i(0)}\right) +\epsilon\\& +
\left(2+\kappa_{10} \ell \right) \sqrt{\frac{\kappa_8}{\rho}} {\tilde\theta_\text{max}} +\kappa_{11} r_0\\
G^{\prime\prime}_{e,3}\triangleq &\epsilon.
\end{split}
\ee
with $\epsilon_1 = \exp(-N)$
\end{thm}
\begin{proof}
The finite time stability result used in \eqref{finitex} still holds for the MMRAC--C. Therefore $G^{\prime\prime}_{x,1}$ in \eqref{GS3} is identical to $G_{x,1}$ in  \eqref{GS}. The Lyapunov function in \eqref{s8:v} has two additional terms in $e_i$ and $\bar \theta$ as compared to the Lyapunov equation in \eqref{can}. Therefore, $G^{\prime\prime}_{e,1}$ now includes the initial conditions of the estimation error $e_i(0)$. $G^{\prime\prime}_{x,2}$ and $G^{\prime\prime}_{e,2}$ are similarly affected. Barbalat Lemma can be used for $G_{e,3}^{\prime\prime}$, and $G_{x,3}^{\prime\prime}$ follows from the same analysis in Appendix \ref{ut3}. The $\eta$ terms arise from the righthand side of the update law in \eqref{s8:law}.
\end{proof}

\section{CMRAC--CO} 
When measurement noise is present, it is often useful to
use a state observer for feedback rather than the plant state.
However, the use of such an observer in adaptive systems
has proved to be quite difficult due to the inapplicability of
the separation principle. In this section, we show how the
CRM can be used to avoid this difficulty for a class of plants.
We denote the resulting adaptive system as CMRAC--CO.

We assume that the plant and reference model dynamics are given by Equations \eqref{s1:sys} and \eqref{s1:ref} with $A_m$ and $L=L_m$ satisfying Equations \eqref{s1:match} and \eqref{s1:ambar}. The control input is now chosen as
\be\label{s6:u}
u= \theta^T(t) x_o + r
\ee
and $x_o$ is the state of the observer dynamics, given by
\be\label{s6:o}
\dot x_o(t) = L_o (x_o(t) - x(t)) + (A_m- b\hat \theta^T(t) )x_o(t) + b u(t).
\ee
Defining $e_m(t) = x(t)-x_m(t)$ and $e_o(t)=x_o-x(t)$, the error dynamics are now given by
\be
\begin{split}\label{s6:e}
\dot e_m(t) =& (A_m+L_m) e_m +b \tilde \theta ^T (t)x_o +b\theta^*e_o \\
\dot e_o(t) =& (A_m+L_o-b\theta^*) e_o -b \bar \theta ^T (t)x_o. 
\end{split}
\ee
For ease of exposition we choose
\be\label{s10:l}
L_m=L_o= L=-\ell  I_{n\times n}.
\ee
The update laws for the adaptive parameters are then defined with the update law
\be\label{s6:law}\begin{split}
\dot \theta = \text{Proj}_\Gamma  (\theta(t),-x_o e_{m}^T P b,f)  - \eta I_{n\times n} \epsilon_\theta\\
\dot {\hat \theta} = \text{Proj}_\Gamma (\hat\theta(t), x_o e_o^T P b,f)+\eta I_{n\times n} \epsilon_\theta
\end{split}\ee
with $\Gamma$ chosen as in \eqref{s1:G}, $\eta>0$, with $P$ from \eqref{s1:lyap}.

\begin{lem} Let 
\be\label{s6:d}
\Delta(\ell) \triangleq    \frac{4 m^2 \norm b \theta_\text{max}^*}{\sigma+2\ell } .\ee
Then, there exists an $\ell^{\prime\prime}$ such that  $0<\Delta(\ell^{\prime\prime})<1$.
\end{lem}

\begin{thm}\label{cmraccosat} Let Assumptions \ref{asm:m} and \ref{asm:thetas} hold with $\ell$ chosen such that $\ell\geq \ell^{\prime\prime}$. Consider the overall CMRAC--CO specified by \eqref{s1:sys}, \eqref{s1:ref}, \eqref{s6:u},  \eqref{s6:o}, \eqref{s6:e} and \eqref{s6:law}. For any initial condition $e_m(0),e_o(0)\in \Re^n$, and $\theta(0)$ and $\hat\theta(0)$ such that $\norm{\theta(0)}\leq\theta_\text{max}$ and $\norm{\hat\theta(0)}\leq\theta_\text{max}$, it can be shown that $e_m(t)$, $e_o(t)$, $\theta(t)$ and $\hat\theta(t)$ are uniformly bounded for all $t\geq0$ and the trajectories in the function 
\be\label{s6:v}
V=   e_m^TP_m e_m + e_o^T P_o e_o + \tilde \theta^T\Gamma^{-1}\tilde\theta + \bar\theta^T \Gamma^{-1} \bar\theta
\ee 
converge exponentially to a set $\mathcal E$ as
\be
\dot V \leq -\alpha_5 V + \alpha_{6}
\ee
where
\be\begin{split}\label{s6:a}
\alpha_5 \triangleq& \frac{\left(1-\Delta(\ell) \right)\left(\sigma + 2 \ell \right)}{m^2}, \\ \alpha_{6} \triangleq&  \frac{2 \left(1-\Delta(\ell ) \right)\left(\sigma + 2 \ell\right)}{\gamma m^2} \tilde\theta_\text{max}^2
\end{split}\ee
and
\ben\begin{split}
\mathcal E \triangleq \left\{ (e_m,e_o,\tilde\theta,\bar\theta) \right | & \norm{e_m}^2  \leq \beta_6 \tilde\theta_\text{max}^2,  \norm{e_o}^2 \leq \beta_6  \tilde\theta_\text{max}^2 \\  & \left.  \norm {\tilde \theta} \leq \tilde\theta_\text{max},\ \norm {\bar \theta} \leq \tilde\theta_\text{max}\right\}
\end{split}\een
with 
\be\label{s6:b}
\beta_6 \triangleq \frac{4(s+l)}{\gamma}.
\ee
\end{thm}
\begin{proof}
see Appendix \ref{pcmraccosat}.
\end{proof}

\subsection{Robustness of CMRAC--CO to Noise }
As mentioned earlier, the benefits of the CMRAC--CO is the use of the observer state $x_o$ rather than the actual plant state $x$. Suppose that the actual plant dynamics is modified from \eqref{s1:sys} as
\be\label{p:n}
\dot x_a(t)= A_p x_a(t)+b u(t), \qquad x(t)=x_a(t)+n(t)
\ee
where $n(t)$ represents measurement noise. For ease of exposition, we assume that $n(t)$ is  bounded and deterministic.

This leads to a set of modified error equations
\be
\begin{split}\label{cne}
\dot e_m(t) =& (A_m+L_m) e_m +b \tilde \theta ^T (t)x_o +b\theta^*e_o  +L_m n(t) \\
\dot e_o(t) =& (A_m+L_o-b\theta^*) e_o -b \bar \theta ^T (t)x_o -L_o n(t)
\end{split}
\ee


\begin{thm}\label{lem:cmraco}
Let Assumptions \ref{asm:m} and \ref{asm:thetas} hold with $\ell$ chosen such that $\ell\geq \ell^{\prime\prime}$. Consider the overall CMRAC--CO specified by \eqref{p:n}, \eqref{s1:ref}, \eqref{s6:u},  \eqref{s6:o}, \eqref{cne} and \eqref{s6:law}. For any initial condition $e_m(0),e_o(0)\in \Re^n$, and $\theta(0)$ and $\hat\theta(0)$ such that $\norm{\theta(0)}\leq\theta_\text{max}$ and $\norm{\hat\theta(0)}\leq\theta_\text{max}$, it can be shown that $e_m(t)$, $e_o(t)$, $\theta(t)$ and $\hat\theta(t)$ are uniformly bounded for all $t\geq0$ and the trajectories in the function $V$ from \eqref{s6:v} converges exponentially 
as
\be
\dot V \leq -\alpha_7 V + \alpha_{8}
\ee
where
\be\begin{split}\label{s11:a}
\alpha_7 \triangleq& \frac{\left(1-\Delta(\ell ) \right)\left(\sigma+ 2\ell \right)}{2m^2} , \\ \alpha_{8} \triangleq&  \frac{\left(1-\Delta(\ell ) \right)\left(\sigma+ 2\ell \right)}{\gamma m^2} \tilde\theta_\text{max}^2\\ &+ \frac{16}{\left(1-\Delta(\ell ) \right)^2} \left(\frac{m^2}{\sigma+2 \ell}\right)^2 \norm{n(t)}^2\end{split}\ee
and\ben\begin{split}
\mathcal E \triangleq \left\{   (e_m,e_o,\tilde\theta,\bar\theta) \right | & \norm{e_m}^2  \leq \beta_6 \tilde\theta_\text{max}^2+\beta_7 \norm{n(t)}^2,  \\  &  \left.\norm{e_o}^2 \leq \beta_6  \tilde\theta_\text{max}^2 +\beta_7 \norm{n(t)}^2, \right.   \\ &  \left.\norm {\tilde \theta} \leq \tilde\theta_\text{max},\ \norm {\bar \theta} \leq \tilde\theta_\text{max} \right\}
\end{split}
\een
with $\beta_6$ defined in \eqref{s6:b} and $\beta_7$ defined as
\be\label{s8:b}
\beta_7 \triangleq \frac{64 m^2 s}{\sigma(1-\Delta(\ell))^3}
\ee
\end{thm}
\begin{proof}see Appendix \ref{applemcmrac}\let\IEEEQED\relax\end{proof}

\subsection{Simulation Study} 
For this study a scalar system in the presence of noise is to be controlled with dynamics as presented in \eqref{p:n}, where $A_p = 1$, $b=1$, and $n(t)$ is a deterministic signal used to represent sensor noise. $n(t)$ is generated from a Gausian distribution with variance 1 and covariance 0.01, deterministically sampled using a fixed seed at 100 Hz, and then passed through a saturation function with upper and lower bounds of 0.1 and -0.1 respectively. The reference model, identifier and observer are from \eqref{s1:ref}, \eqref{s8:o} and \eqref{s6:o} respectively, with $A_m=-1$ and $b=1$.
The controller is defined by \eqref{s6:law}. The design parameters for the two test cases are shown in Table \ref{baus}.
\begin{figure}[t]
\centering
\psfrag{e}[cc][cc][.9]{$e$}
\psfrag{x}[cc][cc][.9]{$x$}
\psfrag{t}[cc][cc][.8]{$t$}
\psfrag{r}[cc][cc][.9]{$x_m$}
\psfrag{openloop}[cl][cl][.8]{open--loop}
\psfrag{closedloop}[cl][cl][.8]{closed--loop}
\psfrag{closedloope0}[cl][cl][.8]{$e^o$: closed--loop}
\psfrag{c1}[cc][cl][.8]{Region 1}

\psfrag{c2}[cc][cl][.8]{Region 2}

\psfrag{c3}[cc][cl][.8]{Region 3}
\includegraphics[width=3.4in]{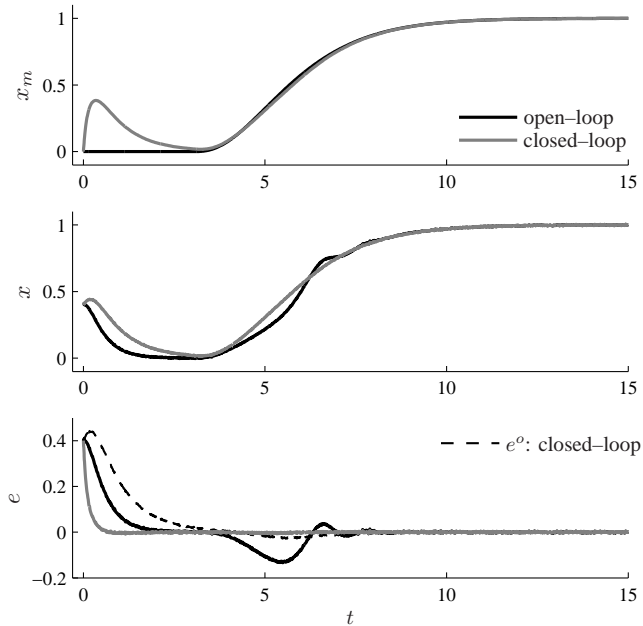}
\caption{(top) reference model trajectories $x_m$, (middle) state $x$, and (bottom)
model following $e$.}\label{fig:ob1}
\end{figure}

\begin{figure}[h]
\centering
\psfrag{u}[cc][cc][.9]{$u$}
\psfrag{u1}[cc][cc][.9]{$\frac{\Delta u}{\Delta t}$}
\psfrag{t}[cc][cc][.8]{$t$}
\psfrag{r}[cc][cc][.8]{$x_m$}
\psfrag{h1}[cc][cc][.8]{$\theta$}

\psfrag{h2}[cc][cc][.8]{$\hat\theta$}
\psfrag{c1}[cc][cl][.8]{Region 1}

\psfrag{c2}[cc][cl][.8]{Region 2}

\psfrag{c3}[cc][cl][.8]{Region 3}

\psfrag{openloop}[cl][cl][.8]{open--loop}
\psfrag{thetastar}[cl][cl][.8]{$\theta^*(t)$}
\psfrag{closedloop}[cl][cl][.8]{closed--loop}
\includegraphics[width=3.4in]{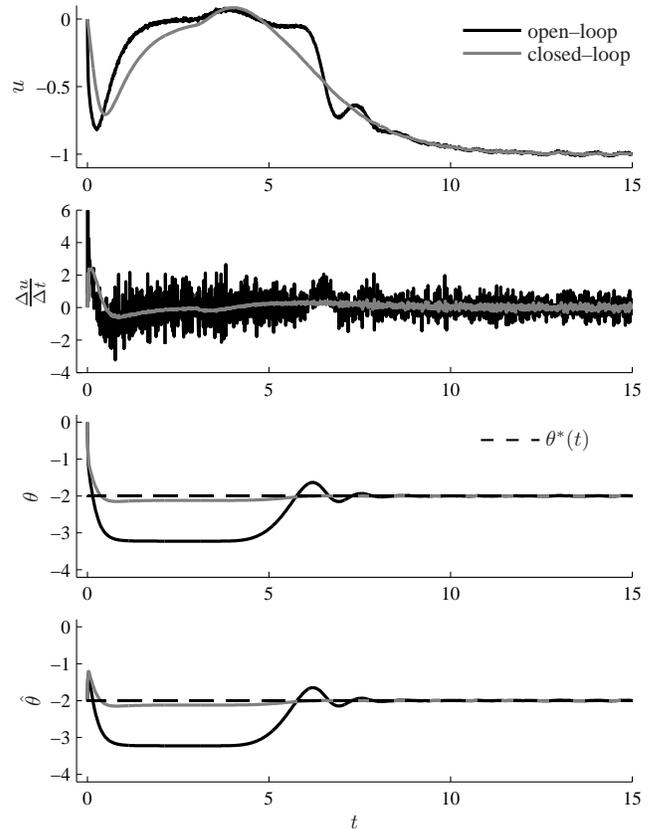}
\caption{(top) Control input $u$, (middle--top) discrete rate of change of control input $\Delta u/ \Delta t$, (middle--bottom) adaptive parameter $\theta(t)$ and (bottom) adaptive parameter $\hat\theta(t)$.}\label{fig:ob2}
\end{figure}

\begin{table}[h]
\caption{Test case free design parameters} \label{baus}
\centering
\begin{tabular}{c c c}
  Paramater & Open--Loop & Closed--Loop  \\ \hline
  $L_{m}$ & 0 & -10 \\
  $L_{i,o}$ & -4 & -4 \\
    $\eta$ & 1 & 1 \\
  $\gamma$ & 100 & 100 \\
   $u(t)$& $\theta x +r$ & $\theta x_o + r $ \\ \hline
\end{tabular}  \end{table}

The simulations have two distinct regions of interest, with Region 1 denoting the first 4 seconds, Region 2 denoting the 4 sec to 15 sec range. In Region 1, the adaptive system is subjected to non--zero initial conditions in the state and the reference input is zero. At $t=4$ sec, the beginning of Region 2, a filtered step input is introduced. Figures \ref{fig:ob1} and \ref{fig:ob2} illustrate the response of the CMRAC--CO adaptive system over 0 to 15 seconds, with $x_m$, $x$, and $e_m$ indicated in Figure \ref{fig:ob1}, and $u$, $\dot u$, $\theta$ and $\hat\theta$ indicated in Figure \ref{fig:ob2}. In both cases, the resulting performance is compared with the classical CMRAC system. The first point that should be noted is a satisfactory behavior in the steady-state of the CMRAC--CO adaptive controller.  We note a significant difference between the responses of CMRAC--CO and CMRAC systems, which pertains to the use of noise free regressors in CMRAC--CO. An examination of $\Delta u/ \Delta t$ in Figure \ref{fig:ob2} clearly illustrates the advantage of CMRAC--CO.

\section{Comments on CMRAC, CMRAC--C and CMRAC--CO}\label{sec:talker}

As discussed in the Introduction, combining indirect and direct adaptive control has always been observed to produce desirable transient response in adaptive control. 
While the above analysis does not directly support the observed transient improvements with CMRAC, we provide a few speculations below: The free design parameter $L_i$ in the identifier is typically chosen to have eigenvalues faster than the plant that is being controlled. Therefore the identification model following error $e_i$ converges rapidly and $\hat\theta(t)$ will have smooth transients. It can be argued that the desirable transient properties of the identifier pass on to the direct component through the tuning law, and in particular $\epsilon_\theta$.

The CMRAC--C differs from classical CMRAC only due to the Luenberger gain $L_m$ in the reference model. Given the contributions of Section \ref{section:perf} which show that the CRM can result in satisfactory transients without the indirect component raises the question if the added complexity of a CMRAC--C is justified. One answer to this question is in the form of the CMRAC--CO, where it is shown that one can design stable observer--based feedback in a CMRAC, allowing noise-free estimation and control.

\section{Conclusion}
This paper concerns the introduction of a feedback gain $L$ in the reference model and the analysis of various adaptive systems with this feature. In particular, we show that with closed-loop reference models (CRM), (i) direct adaptive control structures result in guaranteed transient performance, (ii) combined direct and indirect adaptive controllers result in guaranteed transient performance, and (iii) observer-based feedback can be used in adaptive systems while retaining stability. These are primarily realized using the extra degree of freedom available in the CRM in terms of a feedback gain, and by exploiting exponential convergence properties of the CRM--adaptive system. In all cases, a projection algorithm is used in the adaptive law with a known upper bound on the unknown parameters. 

The main impact of this work is the quantification of transient performance in adaptive systems through $\mathcal L_2$ norms of tracking errors and the control input derivative $\dot u$. It is shown that the introduction of the feedback gain $L$ introduces two time--scales to govern the adaptive system dynamics. The first has to do with the convergence of the tracking error, and the second has to do with adaptation to the unknown plant parameter. By allowing these two time-scales to be separate, transients in the adaptive systems can be controlled without compromising learning of the unknown parameter. This in turn is accomplished by choosing $L$ in an optimal manner. Sub--optimal choices can result in better transients in $e$ only at the expense of slow adaptation leading to a ``water--bed'' effect. This paper, to our knowledge, is the first to illustrate this effect via an exhaustive formal and experimental analysis of CRM--based adaptive systems.

\bibliographystyle{IEEEtran}
\bibliography{ref}

%

\appendices
\section{Projection Operator}\label{ap:proj}
The {\em$\Gamma$--Projection Operator} for two vectors $\theta,y \in \mathbb R^k$, a convex function $f(\theta)\in \Re$ and with symmetric positive definite tuning gain $\Gamma \in \mathbb R^{k \times k}$ is defined as
\begin{equation}\label{eq:proj_vec_gamma}
\text{Proj}_\Gamma(\theta,y,f)=\begin{cases} \Gamma y- \Gamma\frac{\nabla f(\theta) (\nabla f(\theta))^T}{(\nabla f(\theta))^T \Gamma \nabla f(\theta)} \Gamma yf(\theta) \\ \quad \quad \quad \text{ if } f(\theta)>0 \wedge y^T\Gamma \nabla f(\theta)>0\\
\Gamma y  \hspace{8 pt}\quad \quad \text{ otherwise}\end{cases}
\end{equation}
where $\nabla f(\theta)=\left(\frac{\partial f(\theta)}{\partial \theta_1 }\; \cdots \; \frac{\partial f(\theta)}{\partial \theta_k } \right)^T $.  The projection operator was first introduced in \cite{pom92} with extensions in \cite{ioabook} and for a detailed analysis of $\Gamma$--projection see \cite{lav_arXiv}.

\begin{defn}The following compact sets will be referred to in the following analysis:
\be\label{def:d}\begin{split}
\mathcal D_0&\triangleq\{\theta\in\mathbb R^k | f(\theta)\leq 0\}\\
\mathcal D_1&\triangleq\{\theta\in\mathbb R^k | f(\theta)\leq 1\}\\
\mathcal D_\delta&\triangleq\{\theta\in\mathbb R^k | f(\theta)\leq \delta\}.
\end{split}
\ee
\end{defn}

\begin{thm} \label{thm:proj}Given ${\dot \theta= \text{Proj}_\Gamma(\theta,y,f)}$, $f(\theta):\mathbb R^k \rightarrow \mathbb R$ is convex, ${\theta^*\in{\mathcal D}_0}$ and ${\theta(0)\in{\mathcal D}_1}$   
\begin{align}\label{apptheta} \theta(t)&\in{\mathcal D}_1 \forall t\geq 0 \text{ and} \\
\label{eq:prop0vg}
(\theta-\theta^*)^T&(\Gamma^{-1}\text{Proj}_\Gamma({\theta},{y},{f})-y)\leq 0.
\end{align}
\end{thm}
Before we prove the above theorem, we introduce the following two lemmas.

\begin{lem}\label{lem:proj_diff}
Let $f(\theta):\mathbb R^k \rightarrow \mathbb R$ be a continuously differentiable convex function. Choose a constant $\delta>0$. Let $\theta_i$ be an interior point of $\mathcal D_\delta$, defined in \eqref{def:d}. Choose $\theta_b$ as a boundary point so that $f(\theta_b)=\delta$. Then the following holds:
\begin{equation}\label{eq:grad_geo}(\theta_i-\theta_b)^T \nabla f(\theta_b)\leq0\end{equation}
where $\nabla f(\theta_b)=\left(\frac{\partial f(\theta)}{\partial \theta_1 }\; \cdots \; \frac{\partial f(\theta)}{\partial \theta_k } \right)^T $ evaluated at $\theta_b$.
\end{lem}
\begin{proof}
see \cite[Lemma 4]{lav_arXiv}
\end{proof}

\begin{lem}\label{lem:hold10}
Given $\theta^*\in \mathcal D_0$,
\begin{equation}\label{eq:prop1vg}
(\theta-\theta^*)^T(\Gamma^{-1}\text{Proj}_\Gamma({\theta},{y},{f})-y)\leq 0.
\end{equation}
\end{lem}
\begin{proof}
If $ f(\theta)>0 \wedge y^T\Gamma \nabla f(\theta)>0$, then
\begin{equation*}
 (\theta^*-\theta)^T\left(y-\Gamma^{-1}\left(\Gamma y- \Gamma\frac{\nabla f(\theta) (\nabla f(\theta))^T}{(\nabla f(\theta))^T \Gamma \nabla f(\theta)} \Gamma yf(\theta)\right)\right)
\end{equation*}
and using Lemma \ref{lem:proj_diff}
\begin{equation*}
 \frac{(\theta^*-\theta)^T\nabla f(\theta)}{(\nabla f(\theta))^T\Gamma \nabla f(\theta)} {(\nabla f(\theta))^T \Gamma y {f(\theta)} }\leq 0
\end{equation*}
otherwise $\text{Proj}_\Gamma({\theta},{y},{f})=\Gamma y$.
\end{proof}
\begin{proof}[Proof of Theorem \ref{thm:proj}]
We begin by proving \eqref{apptheta}. Consider the function 
\be\label{F}
F(\theta)= f(\theta)^2,\ee and taking its time derivative
\ben
\dot F(\theta) = 2f(\theta)(\nabla f(\theta))^T\dot \theta 
\een
and when $f(\theta)=1$ one has that
\ben
\dot F(\theta) = 2f(\theta) \left(\nabla f(\theta)\right)^T\text{Proj}_\Gamma(\theta,y,f).
\een
With direct substitution of the operator in \eqref{eq:proj_vec_gamma} one finds that \be\label{nablaapp}{(\nabla f(\theta))^T\text{Proj}_\Gamma(\theta,y,f)\leq0}\ee whenever $f(\theta)=1$, and thus \eqref{apptheta} holds. Equation \eqref{eq:prop0vg} is proven with direct application of Lemma \ref{lem:hold10}.
\end{proof}

\begin{thm}\label{thm:projv2}Given 
\be\begin{split}\label{updapp3}
\dot \theta=& \text{Proj}_\Gamma(\theta,y_1,f) - \eta I_{n\times n}(\theta - \hat \theta), \\
\dot {\hat \theta}=& \text{Proj}_\Gamma(\hat\theta,y_2,f) +\eta I_{n\times n} (\theta -\hat\theta)
\end{split}\ee
where $\eta>0$ is a scaler, $\theta^*\in{\mathcal D}_0$, $\theta(0)\in{\mathcal D}_1$, $\hat \theta(0)\in{\mathcal D}_1$ and $f$ is convex
 \be\begin{split} \label{tcb} \theta(t)&\in{\mathcal D}_1 \forall t\geq 0 \text{ and} \\
\hat \theta(t)&\in{\mathcal D}_1 \forall t\geq 0.
\end{split}\ee
\end{thm}
\begin{proof}
Given that $\theta$ and $\hat \theta$ both begin in $\mathcal D_1$ either both parameters hit the boundary of $\mathcal D_1$ simultaneously or only one parameter is at the boundary of $\mathcal D_1$ while the other is strictly inside. Lets consider the case where $\theta(t)$ is on the boundary of $\mathcal D_1$ and thus $f(\theta)=1$ and $\hat \theta \subseteq \mathcal D_1$. Consider the quadratic function $F(\theta)$ as first presented in \eqref{F}. Differentiating $F(\theta)$ and using the update law in \eqref{updapp3} we have
\be\begin{split}\label{Ft2}
\dot F(\theta) =& 2f(\theta) \left(\nabla f(\theta)\right)^T\text{Proj}_\Gamma(\theta,y,f) \\&- \eta 2 f(\theta)\left(\nabla f(\theta)\right)^T(\theta-\hat\theta).
\end{split}\ee
From \eqref{nablaapp} we already know that the first part of \eqref{Ft2} is less than 0. For the second part, given that $f(\theta)$ is convex and since $\hat\theta\in \mathcal D_1$, $\left(\nabla f(\theta)\right)^T(\theta-\hat\theta)\geq 0$, and therefore $\dot F(\theta) \leq 0$. The same result holds for $F(\hat\theta)$, proving \eqref{tcb}.
 \end{proof}

\section{Proof of Lemma \ref{ambarbound}}\label{pf:ambound}

\begin{lem}[{\cite[Lemma 1]{sol94}}]
Any Hurwitz matrix ${A_m\in \Re^{n\times n}}$ with constants $a$ and $\sigma$ as defined in \eqref{s1:ambound} satisfies the following bound for the matrix exponential  
\ben
\norm {\exp(A_m\tau)} \leq m_\delta \exp((-\sigma+\delta a) \tau)
\een
where 
$m_\delta= \frac{3}{2}\left(1+\frac 2 \delta\right)^{n-1}$ and $\delta>0$. The proof follows directly from \cite{sol94}.\end{lem}

\begin{cor}
Setting $\delta = \sigma/(2a)$ the following holds
\be\label{exp2}
\norm {\exp(A_m \tau)} \leq m \exp\left(-\frac \sigma 2 \tau\right),
\ee
where $m =\frac{3}{2} \left(1+4 \varkappa\right)^{n-1}$ and ${\varkappa = \frac a \sigma}$. 
\end{cor}

\begin{lem}
For any diagonal matrix ${L=-lI_{n\times n}}$ the following bound holds for the matrix exponential
\be\label{exp3}
\norm {\exp(L \tau )}  \leq  \exp(- l  \tau) 
\ee
The proof follows from \cite[Section 2]{loa77}.
\end{lem}

\begin{IEEEproof}[Proof of Lemma \ref{ambarbound}(i)]
Beginning with the integral form of Lyapunov's equation in \eqref{s1:lyap}
${P=\int_0^\infty\exp(\bar A_m^T \tau)\exp(\bar A_m \tau)\ d\tau}.$ Due to our choice of $L$, $A_m$ and $L$ commute, thus ${\exp(A_m+L)=\exp(A_m)\exp(L)}$ and
\ben \label{app3:p}
P =   \int_0^\infty \exp(A_m^T\tau) \exp(L^T\tau)  \exp(A_m \tau)\exp(L \tau)\ d\tau.
\een
Using the bound in \eqref{exp2} and \eqref{exp3} the integral just above can be upper bounded and the bound in \eqref{s1:pbound} directly follows.\end{IEEEproof}

\begin{IEEEproof}[Proof of Lemma \ref{ambarbound} (ii)]
Let $\xi\in\Re^n$ be a normalized eigenvector of $P$. By pre-- and postmultiplying \eqref{s1:lyap} by $\xi^T$ and $\xi$, we have 
\ben
\xi^T \bar A_m^TP \xi + \xi^T P\bar A_m \xi =- \xi^TI_{n\times n} \xi
\een
which reduces to 
\ben
 \lambda_i (P) \xi^T (\bar A_m +\bar A_m^T)\xi =-1.
\een
Expanding $\bar A_m$ we have
\ben
 \lambda_i (P) \xi^T ( A_m + A_m^T-2 l I_{n\times n }) \xi =-1.
\een
Finally, using the definitions in \eqref{s1:ambound} and taking the minimum eigenvalue of $P$ we arrive at \eqref{minP} \cite{mor84}.
\end{IEEEproof}

\section{Proof of Theorem \ref{thm:t2}} \label{pf:thm2}
\begin{IEEEproof} Recall the Lyapunov candidate in \eqref{can}, Taking its time derivative one has that
\ben\label{s2:a1}
\dot V \leq  - \norm{e}^2 
 \leq -\frac{1}{\norm P} V + \frac {1} {\norm P \gamma}  \tilde\theta_\text{max}^2.  
\een
Using the upper bound on $P$ from \eqref{s1:pbound} 
\be
\label{a2:a2} 
\dot V \leq  - \alpha_1 V + \alpha_2
\ee
with $\alpha_1$ defined in \eqref{alpha1} and $\alpha_2  \triangleq\frac {\sigma+2\ell} {m^2 \gamma}  \tilde\theta_\text{max}^2$. Using the Gronwall Bellman Inequality, \eqref{a2:a2} implies that 
\be\label{forapp}
V(e,\tilde\theta)\leq \left( V(e(0),\tilde\theta(0))-\frac{\alpha_2}{\alpha_1}\right)\exp(-\alpha_1 t) + \frac{\alpha_2}{\alpha_1}. 
\ee
Thus, $e$ exponentially converges to the set defined by the following inequality
\ben{\lim_{t\rightarrow\infty} {e(t)^TPe(t)} \leq \frac{1} {\gamma} \tilde\theta_\text{max}^2}.\een
Using the bound in Lemma \ref{ambarbound}(ii) we have that 
\be\label{ePe} e^TPe\geq \frac 1 {2(s+\ell) }\norm e^2,\ee 
then we can conclude that
$\lim_{t\rightarrow\infty}\norm {e(t)}^2 \leq \beta_1 \tilde\theta_\text{max}^2 $
where $\beta_1$ is defined in \eqref{b1}. The boundedness of $\theta(t)$ follows from Theorem \ref{thm:proj}.\end{IEEEproof}

\section{Proof of Theorem \ref{ebound}}\label{eboundp}
\begin{proof}
From \eqref{forapp} and \eqref{ePe}, we know that
\ben
\norm {e(t)}^2 \leq  k_0 \exp \left(-\frac{\sigma+2\ell}{m^2} t \right) + k_1
\een
where
\be\label{k01}\begin{split}
k_0 =&   \frac{2(s+\ell) m^2 }{\sigma+2\ell}\norm{e(0)}^2 + \frac{2(s+\ell)}{\gamma } \norm{\tilde\theta(0)}^2  - k_1\\
k_1 = & \frac{2(s+\ell)}{\gamma }\tilde\theta_\text{max}^2 .
\end{split}\ee
Using the following inequalities
\begin{equation*}
 \frac{2(s+\ell) m^2 }{\sigma+2\ell}  \leq \frac{2sm^2 }{\sigma} \text{ and }
 \frac{2(s+\ell)}{\gamma} \leq \frac{2s}{\sigma} \frac{\sigma+\ell}{\gamma}
\end{equation*}
the fact that $\norm{\tilde\theta(0)}\leq\tilde \theta_\text{max}$ and the definition of $\rho$ from \eqref{rhos}, the result in \eqref{e2e} holds with
\be
\kappa_1 =   \frac{2sm^2 }{\sigma} \text{ and }
\kappa_2 =  \frac{2s}{\sigma}.
\ee

Beginning with \be\begin{split}\label{elfirst}
\norm{e(t)}_{L_2}^2 \leq & \int_0^\infty - \dot V(e(t),\tilde\theta(t)) \leq V(e(0),\tilde\theta(0)) \\
\leq& \frac{m^2}{\sigma+ 2\ell} \norm{e(0)}^2 + \frac{1}{\gamma}{\norm{\tilde\theta(0)}}^2,
\end{split}\ee using the definitions of $\rho$ from \eqref{rhos} and the fact that
$
\frac{1}{\sigma+ 2\ell} \leq \frac{1}{\sigma + \ell}
$ the bound in \eqref{el2} holds.
\end{proof}

\section{Proof of Theorem \ref{tbound}}\label{tboundp}
\begin{proof}
Using \eqref{s1:pbound}, the choice for $\Gamma$ in \eqref{s1:G} and the definition of $\rho$ from \eqref{rhos} we have that $
\norm \Gamma \norm P \leq m^2 \rho$. Using the bounds in \eqref{uhm} and \eqref{e2e} for $\norm{x_m(t)}$ and $\norm{e(t)}$ the results in \eqref{t1bound} follow immediately.

For the $\mathcal L_2$ norm we begin by observing that
\ben
\begin{split}\label{bp1}
\norml{\dot\theta(t)}{2}^2 \leq &\norm{\Gamma}^2 \norm{P}^2 \norm b^2\sup\norm {x_m(t)}^2 \int_0^\infty\norm {e(t)}^2  dt \\ &+  \norm{\Gamma}^2 \norm{P}^2 \norm b^2 \sup \norm {e(t)}^2 \int_0^\infty\norm {e(t)}^2 dt.
\end{split}
\een
Taking the supremum of \eqref{uhm} and \eqref{e2e} we have upper bounds for $\sup\norm {x_m(t)}^2$ and $\sup\norm {e(t)}^2$. The $\mathcal L_2$ norm of $e(t)$ is given in \eqref{el2}.\end{proof}

\section{Proof of Theorem \ref{lem:s2}}\label{thm34}
\begin{IEEEproof} Taking the time derivative of the Lyapunov candidate in \eqref{can}, substitution of the update law from \eqref{s2:adaplaw} and the plant dynamics in \eqref{s2:sys}, the derivative of the lyapunov function can be upper bounded as 
\ben
\dot V \leq -  \norm{e}^2 + 2 \norm P \norm  d \norm e+ 2 \frac{ \norm{\dot\theta^*}}{\gamma}\tilde \theta_\text{max}.
\een
After completing the square in $e$ and $d$ we have 
\ben\label{a1}\begin{split}
\dot V \leq& -\frac 1 2  \norm{e}^2 -  \frac 1 2 \left(  \norm e -  2 \norm P \norm d  \right)^2 \\ 
& + 2 \norm P ^2 \norm d^2  + \frac{ 2\normB{\dot\theta^*}}{\gamma}\tilde \theta_\text{max}
\end{split}\een
and then neglecting the negative quantity after $-1/2\norm e^2$,
\ben
 \dot V \leq -\frac 1 2  \norm{e}^2 + 2 \norm P ^2 \norm d^2  + \frac{ 2\norm{\dot\theta^*}}{\gamma}\tilde \theta_\text{max}.
\een
Writing the above inequality in terms of the Lyapunov candidate in \eqref{can} we have
\be\label{a1}
\dot V  \leq -\frac{1}{2 \norm P} V + \frac {1} {2 \norm P \gamma} \tilde\theta_\text{max}^2 + \frac {2\norm{\dot\theta^*}} \gamma \tilde \theta_\text{max}  +2 \norm P ^2 \norm d^2. 
\ee
Using the upper bound on $\norm P$ from \eqref{s1:pbound}  and rewriting \eqref{a1} in terms of the design parameters $\gamma$ and $\ell$ we have
\be
\dot V \leq -\alpha_3 V + \alpha_4
\ee
where $\alpha_3$ is defined just before $\alpha_4$ in \eqref{alpha4}.
Following the same procedure as in Appendix \ref{pf:thm2} we conclude that
\be\label{ePed}
\lim_{t\rightarrow\infty}{e^TPe} \leq \frac 1 \gamma \tilde\theta^2_\text{max} +  \frac {4m^2\norm{\dot\theta^*}} {(\sigma+2\ell)\gamma} \tilde \theta_\text{max}  +2 \left(\frac{m^2}{\sigma+2\ell}\right) ^3 \norm d^2.\ee
Recalling the fact that ${s\geq\sigma>0}$ from Lemma \ref{sgs} we can conclude that
\be\label{spl}
\frac{2(s+\ell)}{\sigma+2\ell}\leq \frac{2s}{\sigma}.
\ee
Using the bound above along with that in \eqref{ePe} the inequality in \eqref{ePed} can be simplified as 
\be
\lim_{t\rightarrow\infty}\norm {e(t)}^2 \leq \beta_1 \tilde\theta^2_\text{max} + \beta_2 \norm{\dot\theta^*}  \tilde\theta_\text{max}+\beta_3 \norm d^2
\ee
where $\beta_1$ is defined in \eqref{b1}, and $\beta_2$ and $\beta_3$ are defined in \eqref{b2}. The boundedness of $\theta(t)$ follows from Theorem \ref{thm:proj}.\end{IEEEproof}

\section{Proof of Theorem \ref{xmbound}}\label{xmboundp}
\begin{IEEEproof}
The dynamics of the CRM and the ORM are given in \eqref{s1:ref} and \eqref{eo} respectively and leed to the following
\be
\dot x_m(t) - \dot x_m^o(t) = A_m (x_m(t)-x_m^o(t)) - L e.
\ee
Given that the reference model will have the same initial condition regardless of being closed or open, we then have that 
\be
\norm{x_m(t)-x_m^o(t)}  \leq m \int_0^t\exp({-\frac{\sigma}{2}(t-\tau)}) \ell e(\tau) d\tau 
\ee
where the matrix exponential bound came from \eqref{exp2}. Using the Cauchy--Schwartz inequality we have the following bound
\be
\norm{x_m(t)-x_m^o(t)}  \leq \frac{\ell m}{\sqrt\sigma} \norml{e(t)}{2}.
\ee \let\IEEEQED\relax
\end{IEEEproof}

\section{Proof of Theorem \ref{thmudot}}\label{tan}
Taking the time derivative of $u$ in \eqref{s1:u}
\be \begin{split}
\dot u(t)=& - b^T P e(t) x^T(t)\gamma I_{n\times n} x(t)\\ &+ \theta^T\left(A_m x(t)  + b \left( \tilde\theta^T  x(t)+ r(t)\right)\right)+\dot r(t).
\end{split}\ee
Substitution of the upper bound on $P$ from \eqref{s1:pbound}, using the definition of $a_\theta$ from \eqref{es} and the bounds on the reference trajectory from Assumption \ref{asst1} results in the following bound
\be\begin{split}\label{s1:udotbound}
\abs{\dot u(t)} \leq &\frac{m^2\gamma}{\sigma+2\ell}\norm b \norm{e(t)} \norm{x(t)}^2  \\ & + \theta_\text{max}\left(a_\theta \norm {x(t)}+r_0\right) + r_1.
\end{split}\ee
\subsection{Proof of Theorem \ref{thmudot}, $t\in\mathbb T_1$}
The following Lemma is useful:
\begin{lem}\label{lemfinx}[Finite time stability]
If $r$ satisfies Assumption \ref{asst1}, then
\be\label{finitex}
\norm{x(t)}  \leq  \norm{e(0)} \exp \left( a_\theta t \right) + \frac{\norm b r_0}{a_\theta}\left ( \exp(a_\theta t) - 1 \right ),\quad t\geq0
\ee
where $a_\theta$ is defined in \eqref{es}.
\end{lem}
\begin{proof}
Suppose $z(t)\in \Re$ is defined as the solution to
\be\label{zd}
\dot z(t)=a_\theta z(t) + \norm b r_0.
\ee
It can be shown that if $z(0)=\norm{x(0)}$, then
\be\label{zb}
\norm{x(t)}\leq z(t)\quad\forall\ t\geq0
\ee using \cite[Theorem 8.14]{dan70}.\end{proof} Using the fact that $x(0)=e(0)$ which follows from Assumption \ref{assxm0}, Lemma \ref{lemfinx} and the definitions of $a_\theta$ and $\tau_1$ we obtain that
\be\label{xt1}
\sup_{t\in \mathbb T_1} \norm{x(t)} \leq G_{x,1}
\ee
where $G_{x,1}$ is defined in \eqref{GS}.

Beginning with \eqref{e2e}, taking the square root of the expression and noting that ${\sqrt{c_1+c_2}\leq \sqrt{c_1}+\sqrt{c_2}}$ for all ${c_1,c_2>0}$, we obtain
\be\label{egron}
\norm{e(t)} \leq \sqrt{\kappa_1} \exp \left(-\tfrac{1}{\tau_1} t \right) \norm{e(0)} + \sqrt{\frac{\kappa_2}{\rho}} {\tilde\theta_\text{max}}
\ee
where $\tau_1$ is defined in \eqref{taudef1}. This verifies that
\be\label{e4t}
\sup_{t\in\mathbb T_1} \norm{e(t)} \leq G_{e,1}
\ee 
where $G_{e,1}$ is defined in \eqref{GS} . Using \eqref{s1:udotbound}, \eqref{xt1}, and \eqref{e4t}, Theorem \ref{thmudot} for $t\in\mathbb T_1$ is proved.

\subsection{Proof of Theorem \ref{thmudot}, $t\in\mathbb T_2$}
From \eqref{egron}  it is easy to see that,
\be\label{e4t2}
\sup_{t>N\tau_1} \norm{e(t)} \leq G_{e,2}
\ee 
where $G_{e,2}$ is defined in \eqref{GS}.

From \eqref{s1:ref} and the bound on $\exp(A_m t)$ in \eqref{exp2}, we have that
\be\label{axm10}
\norm{x_m(t)}\leq  m \int_0^t  \exp\left( -\tfrac{1}{\tau_2} (t-\tau) \right) \left( l \norm{ e(\tau)}+\norm b \norm r\right) d \tau \ee 
Using the integral transform of LTI systems, the bound for $\exp(A_m)$ from \eqref{exp2}, the bound for $\norm{e(t)}$ from \eqref{egron}, \eqref{axm10} takes the form
\be\label{uhm2}\begin{split}
\norm {x_m(t)} \leq & m_1 \norm {e(0) }\left( \exp\left(-\tfrac{1}{\tau_2}t\right) -\exp\left(-\tfrac{1}{\tau_1}t\right) \right) \\ &+ \frac{ 2lm }{\sigma} \sqrt{\frac{2(s+\ell)}{\gamma}}\tilde\theta_\text{max}\left( 1-\exp\left(-\tfrac{1}{\tau_1}t\right) \right) \\
&+ \frac{ 2 \norm{b} m} {\sigma} r_0  \left( 1-\exp\left(-\tfrac{1}{\tau_1}t\right) \right)
\end{split}\ee
where $m_1 \triangleq \frac{2lm^4\sqrt{\frac{2s}{\sigma}}}{\sigma+2\ell-\sigma m^2}$. 
 
Given that $x=e+x_m$, using \eqref{e4t} and \eqref{uhm2} one can conclude that 
\be\label{tt12}
\sup_{t\geq N\tau_1} \norm{x(t)}\leq G_{x,2}
\ee
where $G_{x,2}$ is defined in \eqref{GS}. Using \eqref{s1:udotbound}, \eqref{e4t2}, and \eqref{tt12}, Theorem \ref{thmudot} for $t\in\mathbb T_2$ is proved.

\subsection{Proof of Theorem \ref{thmudot}, $t\in\mathbb T_3$} \label{ut3}
$G_{e,3}$ follows from Corollary \ref{Tdef}. $G_{x,3}$ follows from \eqref{uhm2}, where it is noted that $t\geq N\tau_2$, and the fact that $\norm x \leq \norm e + \norm{x_m}$
\section{Proof of Theorem \ref{propl0}}\label{pfl0}
The bound given for $x(t)$ and $e(t)$ over the time period $[0,N\tau_1^*]$ in \eqref{finitex} holds regardless of the choice of $\ell$. Thus the bound in \eqref{ub1} holds for $\ell=0$, and therefore the bound in \eqref{uox2} for $t\in[0,N\tau_1^*]$ is the same as that in \eqref{uox} where $\rho$  has been replaced by $\rho_0$. The Gronwall--Bellman analysis used to obtain the bound for $e(t)$ would follow with a similar bound to that in \eqref{egron} where the exponent would now have the time constant $\tau_2^\prime$ with $e(t)$ exponentially decaying to $\sqrt{1/\rho_0}\tilde\theta_\text{max}$. For $t>N\tau_2^\prime$, $\norm{e(t)}$ would have decayed past 4 time constants. Therefore, the coefficient $G_{e,2}$ would apply for the ORM case when $t>N\tau_2^\prime$ and $\ell=0$. The bound for $x(t)$ would not contain the parameter $\ell$. Therefore, \eqref{uox2} for $t> N\tau_2^\prime$ is identical in structure to \eqref{uox} for $t>N\tau_1$ with $\ell=0$ and $\rho$ being replaced with $\rho_0$. The asymptotic properties of the adaptive system hold regardless of the choice of $\ell$ and therefore Corollary \ref{Tdef} holds when $l=0$ as well and thus the bounds in \eqref{uox2} for $t\geq T_1^\prime$ hold as well.

\section{Proof of Theorem \ref{cmracsat}} \label{proofcmracsat}
\begin{proof}
Taking the time derivative of $V$ in \eqref{s8:v} results in
\be
\dot V \leq -\norm{e_m}^2 -\norm{e_i}^2 -2 \frac\eta\gamma\epsilon_\theta^2.
\ee
Substitution of $V$ in \eqref{s8:v} results in
the bound in \eqref{tbone3}. Using the bound in Lemma \ref{ambarbound}--(ii) we have that 
\ben e_m^TP_me_m \geq \frac 1 {2(s+\ell) }\norm {e_m}^2 \text{ and } e_i^TP_ie_i \geq \frac 1 {2 (\sigma+\ell)  }\norm {e_i}^2 \een 
then we can conclude that
${\lim_{t\rightarrow\infty}\norm {e_m(t)}^2 \leq \beta_4 \tilde\theta_\text{max}^2}$ and ${\lim_{t\rightarrow\infty}\norm {e_i(t)}^2 \leq \beta_5 \tilde\theta_\text{max}^2}$
where $\beta_4$ and $\beta_5$ are defined in \eqref{s4:b}. The boundedness of $\theta(t)$ and $\hat\theta(t)$ follows from Theorem \ref{thm:projv2}. The asymptotic limit to zero comes from the application of Barbalat Lemma.
\end{proof}

\section{Proof of Theorem \label{eperf}}\label{proofeperf}
The bounds in \eqref{e2me} and \eqref{e2i} follow from the application of Gronwall--Bellman to the result in \eqref{tbone3} with the lower bound for $\min \lambda_i(P)$ in \eqref{minP} and the change of parameters from \eqref{rhos} being used.

Beginning with \be\begin{split}\label{elfirst}
\norm{e_m(t)}_{L_2}^2 \leq & \int_0^\infty - \dot V(e(t),\tilde\theta(t)) \leq V(e(0),\tilde\theta(0)) \\
\leq& \frac{m^2}{\sigma+ 2\ell} \norm{e_m(0)}^2 +\frac{1}{2(\sigma+\ell)} \norm{e_i(0)}^2\\&+ \frac{2}{\gamma}{\norm{\tilde\theta(0)}}^2,
\end{split}\ee using the definitions of $\rho$ from \eqref{rhos}, the fact that
$\frac{1}{\sigma+ 2\ell} \leq \frac{1}{\sigma + \ell}$ the bound in \eqref{eml2} holds. This same approach can be used to obtain the bound in \eqref{eei2}.

\section{Proof of Theorem \ref{cmraccosat}}\label{pcmraccosat}
\begin{proof}
Taking the time derivative of $V$ in \eqref{s6:v} results in
\be
\dot V \leq -\left(1-\Delta(\ell ) \right)\left(\norm{e_m}^2 +\norm{e_o}^2\right) -2\frac\eta\gamma\epsilon_\theta^2.
\ee
where $\Delta(l)$ is defined in \eqref{s6:d}. Substitution of $V$ in \eqref{s6:v} results in
\be
\dot V \leq -\alpha_5 V + \alpha_{6}
\ee
where $\alpha_5$ and $\alpha_{6}$ are defined in \eqref{s6:a}. Using the bound in Lemma \ref{ambarbound}--(ii) we have that 
\ben e_m^TPe_m \geq \frac 1 {2(s+\ell) }\norm {e_m}^2 \text{ and } e_o^TPe_o \geq \frac 1 {2(s+\ell) }\norm {e_o}^2 \een 
then we can conclude that
${\lim_{t\rightarrow\infty}\norm {e_m(t)}^2 \leq \beta_6 \tilde\theta_\text{max}^2}$ and ${\lim_{t\rightarrow\infty}\norm {e_o(t)}^2 \leq \beta_6 \tilde\theta_\text{max}^2}$
where $\beta_6$ is defined in \eqref{s6:b}. The boundedness of $\theta(t)$ and $\hat\theta(t)$ follows from Theorem \ref{thm:projv2}.
\end{proof}

\section{Proof of Theorem \ref{lem:cmraco}}\label{applemcmrac}
\begin{proof}
Taking the time derivative of $V$ in \eqref{s6:v} results in \be
\begin{split}\dot V \leq& -\left(1-\Delta(\ell ) \right)\left(\norm{e_m}^2+\norm{e_o}^2\right) -2\frac\eta\gamma\epsilon_\theta^2 \\&+2 \norm P \norm{n(t)}\norm{e_m(t)} +2 \norm P \norm{n(t)}\norm{e_o(t)}
\end{split}.\ee
completing the square in $\norm{e_m}\norm n$ and $\norm{e_o} \norm{n}$
\be\begin{split}
\dot V \leq &-\frac{\left(1-\Delta(\ell ) \right)}{2} \left(\norm{e_m}^2 +\norm{e_o}^2\right) -2\frac\eta\gamma\epsilon_\theta^2 \\ 
&-\frac{\left(1-\Delta(\ell ) \right)}{2} \left(\norm{e_m} -  \frac{4} {\left(1-\Delta(\ell ) \right)} \norm P \norm{n(t)}\right)^2 \\  
&-\frac{\left(1-\Delta(\ell ) \right)}{2} \left(\norm{e_o} -  \frac{4} {\left(1-\Delta(\ell ) \right)} \norm P \norm{n(t)}\right)^2 \\
& + \frac{16}{\left(1-\Delta(\ell ) \right)^2} \norm P^2 \norm{n(t)}^2 .
\end{split}
\ee
Neglecting the negative terms in lines 2 and 3 from the equation above and substitution of the norm for $P$ we have that
\be\begin{split}
\dot V \leq &-\frac{\left(1-\Delta(\ell ) \right)}{2} \left(\norm{e_m}^2 +\norm{e_o}^2\right) -2\frac\eta\gamma\epsilon_\theta^2 \\ 
& + \frac{16}{\left(1-\Delta(\ell ) \right)^2} \norm P^2 \norm{n(t)}^2 .
\end{split}
\ee
which in terms of $V$ is identical to
\be\begin{split}
\dot V \leq &-\frac{\left(1-\Delta(\ell ) \right)\left(\sigma+ 2\ell \right)}{2m^2}  V + \frac{\left(1-\Delta(\ell ) \right)\left(\sigma+ 2\ell \right)}{\gamma m^2} \tilde\theta_\text{max}^2 \\ 
& + \frac{16}{\left(1-\Delta(\ell ) \right)^2} \left(\frac{m^2}{\sigma+2 \ell}\right)^2 \norm{n(t)}^2 .
\end{split}\ee
\be
\dot V \leq -\alpha_7 V + \alpha_{8}
\ee
where $\alpha_7$ and $\alpha_{8}$ are defined in \eqref{s11:a}. Using the bound in Lemma \ref{ambarbound}--(ii) we can conclude that
${\lim_{t\rightarrow\infty}\norm {e_m(t)}^2 \leq \beta_6 \tilde\theta_\text{max}^2+\beta_7 \norm{n(t)}^2}$ and ${\lim_{t\rightarrow\infty}\norm {e_o(t)}^2 \leq \beta_6 \tilde\theta_\text{max}^2+\beta_7 \norm{n(t)}^2}$
where $\beta_7$ is defined in \eqref{s8:b}. The boundedness of $\theta(t)$ and $\hat\theta(t)$ follows from Theorem \ref{thm:projv2}.
\end{proof}

\end{document}